\documentclass[a4paper,12pt]{article}
\usepackage{amsmath,amssymb}
\usepackage{amsthm}
\usepackage{graphicx,tikz}
\usepackage{enumerate}
\usepackage{titlesec}
\usepackage{epstopdf}
\usepackage{caption}
\usepackage{xcolor}
\usepackage{accents}

\usepackage{wallpaper}

\newtheorem{Thm}{Theorem}[section]
\newtheorem*{Thm*}{Theorem}
\newtheorem{Prop}[Thm]{Proposition}
\newtheorem{Lem}[Thm]{Lemma}
 
\newtheorem{Cor}[Thm]{Corollary}
\newtheorem{Fact}[Thm]{Fact}

\newtheorem{Exa}[Thm]{Example}

\theoremstyle{remark}

\theoremstyle{definition}

\newtheorem{Def}[Thm]{Definition}

\newtheorem*{Def*}{Definition}

\numberwithin{equation}{section}

\newcommand{\m}[1]{\mathbb{ #1}}
\newcommand{\mc}[1]{\mathcal{ #1}}

     \def\ol{\overline}    
   
\def\al{\alpha}       \def\be{\beta}        
\def\de{\delta}       \def\eps{\varepsilon}  
       \def\la{\lambda}      
\def\si{\sigma}                
\def\ph{\varphi}               
              \def\De{\Delta}
      \def\Si{\Sigma}       \def\Ph{\Phi}
                  \def\Om{\Omega}

\def\bB{\partial B}

\def\bOm{\partial \Om}

\def\bX{\partial X}

\def\GGG{{\bf GGG }}

\theoremstyle{definition}

\theoremstyle{remark}

\newtheorem{Rmq}[Thm]{Remark}

\numberwithin{equation}{section}
\newfont{\goth}{eufm10 at 12pt}
\newfont{\gots}{eufm8 at 9pt}

\def\bt{\begin{Thm}}
\def\et{\end{Thm}}
\def\br{\begin{Rmq}}
\def\er{\end{Rmq}}

\def\bc{\begin{Cor}}
\def\ec{\end{Cor}}
\def\bp{\begin{Prop}}
\def\ep{\end{Prop}}
\def\bl{\begin{Lem}}
\def\el{\end{Lem}}
\def\bd{\begin{Def}}
\def\ed{\end{Def}}
\def\bq{\begin{quotation}}
\def\eq{\end{quotation}}
\def\bfa{\begin{Fact}}
\def\efa{\end{Fact}}
\def\bexa{\begin{Exa}}
\def\eexa{\end{Exa}}

\def\ra{\rightarrow}
\def\vs{\vspace{1em}}

\begin{document}
\title{Bounded harmonic maps}
\author{Yves Benoist \&
Dominique Hulin}
\date{}

\maketitle

\begin{abstract} 
The classical Fatou theorem identifies
bounded harmonic functions on the unit disk  with
bounded measurable functions on the boundary circle. 
We extend this theorem to bounded harmonic maps. 
\end{abstract}

\renewcommand{\thefootnote}{\fnsymbol{footnote}} 
\footnotetext{     
\renewcommand{\thefootnote}{\arabic{footnote}} }

\newpage 

{\footnotesize \tableofcontents}
\newpage



\section{Introduction}

\bq 
The aim of this paper is to present 
an extension to harmonic maps 
of a  classical theorem for harmonic functions 
due to Fatou around 1905-1910.
\eq

\subsection{The Fatou theorem}
\label{secfatouherglotz}
\bq 
We first recall the classical theorem of Fatou for bounded harmonic functions on the Euclidean disk. 
\eq 

This theorem deals with the unit open ball $B$ and with the unit sphere $S=\bB$ in the Euclidean space $\m R^k$ with $k=2$ or, more generally, with $k\geq 2$.
It identifies the  space $\mc H_b(B,\m R)$ of bounded harmonic functions $h:B\ra \m R$ on the ball
with the space $L^\infty(\bB,\m R)$ of bounded measurable functions on the boundary $\bB$.
We recall that a harmonic function on $B$ is a $C^2$-function $h$ that satisfies $\Delta_0 h=0$, where $\Delta_0$ is the Euclidean Laplacian. We denote by $\si_0$  the rotationally invariant
probability measure on the sphere $\bB$, and we refer to Section \ref{secmaindefinition} for the definition of a non-tangential limit.

\bfa 
\label{facfatou}
{\bf (Fatou)} 
$a)$ Let $h:B\to\m R$ be a bounded harmonic function. \\ For $\si_0$-almost all $\xi$ in $\bB$, the function $h$ admits a non-tangential limit $\ph(\xi):=\underset{x\ra\xi}{\rm NTlim}\,h(x)$ at the point $\xi$.\\
$b)$ The map $h\mapsto \ph$  is a bijection $\beta : \mc H_b(B,\m R)\ra L^\infty(\bB,\m R)$ called the boundary transform.
\efa

The inverse of the map $\beta$ is given by an explicit formula, the Poisson formula.
For every $\ph\in L^\infty(\bB,\m R)$, one can indeed recover $h$ as  $h=P_0\ph$ where 
$P_0\ph$ is the bounded harmonic function defined on $B$  by
$$
P_{0}\ph(x):=\int_{\bB}  \ph(\xi)\,P_{0,\xi}(x)\,{\rm d}\si_0(\xi),
\;\;{\rm where}\;\;
P_{0,\xi}(x)=\frac{1-|x|^2}{|x-\xi|^k }$$ 
is  the Poisson kernel. 

The proof of this fact
can be found in  Rudin's book \cite[Chap. 11]{RudinRealComplex} or in \cite{Koosis} when $k=2$,
or in Armitage and Gardiner's book \cite[Chap. 4]{ArmitageGardinerPotential}
for $k\geq 2$.

Proving extensions of this fact has
a long history that has already lasted for more than a century. 
Indeed, an important  goal of  Potential Theory is to understand to what extent
this fact still holds   for either harmonic or superharmonic functions,  on more general spaces.

\vspace{0.5em}
The aim of this paper is to extend Fatou's theorem to bounded harmonic maps. 
We will allow the target space $Y$ to be  
any complete CAT(0) space,  the first examples to have in mind being the hyperbolic spaces $\m H^k$.

We will also allow more general source spaces.
Since, in dimension $k=2$, the harmonicity condition
depends only on the conformal structure on the source space $B$, 
we can think of $B$ as the hyperbolic plane.
We will explain in  Theorem \ref{thmbijectiveboundary} how to replace $B$ by a 
\GGG Riemannian manifold $X$.
Later on, we will also explain  in Corollary \ref{corbijectiveboundary} 
how to 
replace $B$
by any bounded  Riemannian domain
$\Om$ with Lipschitz boundary.

\subsection{Main result}
\label{secmainresult}

\bq 
We now state  our main result, postponing the definitions to Section \ref{secmaindefinition}.
\eq

\bd
\label{defGGG}
We will say that a Riemannian  manifold $X$ is  \GGG as a shortcut for   
{\bf G}romov Hyperbolic
with Bounded {\bf G}eometry and 
Spectral {\bf G}ap.
\ed

\bt
\label{thmbijectiveboundary}
Let $X$ be a  \GGG Riemannian manifold,\! and $Y$ be a proper {\rm CAT(0)}-space.\\
$a)$ Let  $h\!:\! X\!\ra\! Y$ be a bounded harmonic map. 
Then, for $\si$-almost all $\xi\!\in\! \bX$,
the map $h$ admits a non-tangential limit $\ph(\xi):=\underset{x\ra\xi}{\rm NTlim}\,h(x)$ at the point $\xi$.\\
$b)$ The map $h\mapsto \ph$
is a bijection 
$
\beta:\mc H_b( X,Y){\ra} L^\infty(\bX,Y).
$
\et

The main examples of \GGG Riemannian manifolds $X$ are the pinched Hadamard manifolds ~: those have negative curvature.
The condition \GGG allows a little bit of positive curvature on $X$. It also allows $X$ to be non-contractible. For instance, the quotient of 
a pinched Hadamard manifold by a convex cocompact 
group of isometries is \GGG.

\subsection{Main definitions}
\label{secmaindefinition}

Here are the definitions that are needed to understand our theorem \ref{thmbijectiveboundary}.
All our manifolds will be  assumed to be connected and
with dimension $k\geq 2$.

\bd
\label{defboundedgeometry} 
A Riemannian manifold $X$ has {\it bounded geometry}
if it is complete, with bounded sectional curvature   $-K_{max}\leq K_X\leq K_{max}$, and if the 
injectivity radius has a uniform lower bound
${\rm inj}_X\geq r_{min}>0$.
\ed

As explained in \cite[Section 1.1]{KemperLohkamp}, one could replace 
in Definition \ref{defboundedgeometry} the bound on the sectional curvature
by a bound on the Ricci curvature. Indeed, the important features of bounded geometry also hold if we just have a bound on the Ricci curvature.

\bd
\label{defgromov} 
The Riemannian manifold $X$ is Gromov hyperbolic
if there exists $\de>0$ such that, for all
$o$, $x$, $y$, $z$ in $X$ one has 
\begin{equation}
\label{eqngromovproduct}
(x|z)_o\geq \min ( (x|y)_o,(y|z)_o )-\de.
\end{equation}
Here $(x|y)_o:=\tfrac12(d(o,x)+d(o,y)-d(x,y))$
is the Gromov product of the points $x$ and $y$ seen from $o$.

In this case $\bX$ will denote the ``Gromov boundary'' of $X$ and $\ol X=X\cup\bX $ will be the ``Gromov compactification'' of $X$ where $\bX$ is the set of geodesic rays on $X$, two geodesic rays being identified 
if they remain within bounded distance from each other.
See \cite{GhysHarp90}. 
\ed

\bd
\label{defspectralgap}
A Riemannian manifold $X$ has a {\it spectral gap},
or  a {\it coercive Laplacian},
if the Rayleigh quotients admit a uniform lower bound
\begin{equation}
\la_1:=\inf\limits_{\ph\in C^\infty_c(X)}\frac{\int_X\|\nabla\ph\|^2\,{\rm d}v_g}{\int_X\ph^2\,{\rm d}v_g}
\; >\; 0.
\end{equation}
\ed
Note that  the spectral gap implies 
that $X$ is non-compact.

\bd 
\label{defHadamard}
A pinched Hadamard manifold $X$ is
a complete simply-connected Riemannian manifold with dimension at least $2$
whose sectional curvature is pinched between two negative constants~: 
$-b^2\leq K_X\leq -a^2<0$.
\ed
Examples are~: hyperbolic spaces $\m H^k$, 
rank one non-compact Riemannian symmetric spaces,
any small perturbation of those...

\bd 
A CAT(0) space $Y$ is a geodesic metric space such that,
for every geodesic triangle $T$ in $Y$, 
 there exists a $1$-Lipschitz map $j:T_0\to T$ where $T_0$ is the
triangle  the Euclidean plane with same side lengths as $T$ and  $j$
 sends each vertex of $T_0$ to the corresponding vertex of $T$.
See \cite{BridsonHaefliger}.
\ed

Examples are~: Hadamard manifolds (namely, complete and simply connected Riemannian manifolds with non positive curvature),  Euclidean buidings,
$\m R$-trees, convex subsets in Hilbert spaces... 

It is not restrictive to assume that $Y$ is complete, since the metric completion of a CAT(0) space still is a CAT(0) space. 
Since the closed balls $B(y_0,R)$ in $Y$ are also CAT(0) spaces, 
this will allow us to assume in the proof of 
Theorem \ref{thmbijectiveboundary} that $Y$ is bounded.
We will sometimes assume that $Y$ is proper, i.e. 
that these closed balls $B(y_0,R)$ are compact.

\bd 
A map $h: X\ra Y$ is (energy minimizing) harmonic if it is locally Lipschitz continuous 
and
if it is a minimum for the Korevaar-Schoen energy $E(h)$ with respect to variations of $h$ with compact support $Z\subset  X$.
\ed

When $Y$ is a CAT$(0)$ Riemannian manifold, the Korevaar-Schoen energy on $Z$ coincides with the
Dirichlet energy  $E(h)=\int_Z |Dh(x)|^2 dv_{g}(x)$. 
In this case, the harmonicity condition can be expressed by a partial differential equation which  is not linear  any more,
see \cite{EellsSampson64}, \cite{Hamilton75} or \cite{Jost84}. 
When $Y$ is only a {\rm CAT}(0) space, the energy of $h$ on $Z$ is the integral  $E(h)=\int_Z e_h(x) dv_{g}(x)$ of the energy density $e_h$ where, 
for a Lipschitz continuous map $h$, the energy density is  given by
$e_h(x)=\limsup\limits_{\eps\ra 0}\eps^{-2-k}v_k^{-1}\int_{B(x,\eps)} d(h(x),h(x'))^2dv_{g}(x') $, where $v_k$ is the volume of the 
unit Euclidean ball and
where this limit should be understood in a weak sense.
See \cite[Section 1.5]{KorevaarSchoen1} for a precise definition.
See also \cite{Jost95}.\vs 

The measure $\si$ refers to any finite Borel measure on $\bX$ which is equivalent to 
the harmonic measures on $\bX$.
The ``$\si$-almost surely'' means that the property holds except 
on a set of measure zero for the  measure $\si$ on $\bX$.
Note that, when $X$ is a pinched Hadamard manifold, such a measure $\si$ is often found to be singular with respect to  the Lebesgue measure on the sphere $\bX$.

\vs
The set
$\mc H_b( X,Y)$ is the set of bounded
harmonic maps $h: X\ra Y$, and the set
$L^\infty(\bX,Y)$ is the set of bounded measurable maps from $\bX$ to $Y$ where two measurable maps are identified  if they are $\si$-almost surely equal.

\bd 
A function $h: X\ra Y$ has a non-tangential limit 
$y$ at a point $\xi\in \bX$ 
(also called a conical limit), and we write 
$y=\underset{x\ra\xi}{\rm NTlim}\,h(x),$
if $\displaystyle y=\lim_{n\ra\infty} h(x_n)$ holds for any sequence $(x_n)$ in $ X$ 
converging non-tangentially to $\xi$, i.e. such that
$\;\displaystyle\sup\limits_{n\geq 1}d(x_n, o\xi)<\infty$
where $o\xi$ is  any geodesic ray from a point $o\in X$ to $\xi$.
\ed

\subsection{Previous results}
\label{secprevious}

When $Y=\m R$, we are dealing with harmonic functions.  
As we have already seen in Section \ref{secfatouherglotz},
Theorem \ref{thmbijectiveboundary} for $X=B$ is  the classical 
Fatou  theorem.  
The extension to the case
where  $X$ is a pinched Hadamard manifold appeared in the $80's$
and is due to Anderson and Schoen in  \cite{AndersonSchoen}.
The extension to the case where $X$ is a \GGG Riemannian manifold
is due to Ancona in \cite{AnconaStFlour}.

When $Y$ is a CAT$(0)$ Riemannian manifold 
and $X$ is a pinched Hadamard manifold, Theorem \ref{thmbijectiveboundary}.$a$ is due to
Aviles, Choi, Micallef in \cite[Thm 5.1]{AvilesChoiMicallef91},
and Theorem \ref{thmbijectiveboundary}.$b$ is conjectured to be true 
by these authors. Indeed, 
as a final observation 
in \cite[Section 1]{AvilesChoiMicallef91} they write that 
such a theorem would  be ``a consequence of the solvability 
of the Dirichlet problem with $L^\infty$ boundary condition''. This solvability 
is one of the main technical issues in our paper (Proposition \ref{prosurjectiveboundary}).
Note that the solvability of the Dirichlet problem with continuous boundary 
condition is proven in \cite[Thm 3.2]{AvilesChoiMicallef91}. 
The first case of Theorem  \ref{thmbijectiveboundary}.$b$ 
that seems to be new is when both
$X$ and $Y$ are the hyperbolic plane $\m H^2$.

When $Y$ is a  CAT(0) space, the proof of Theorem \ref{thmbijectiveboundary} will rely on the  solution of the Dirichlet problem 
for harmonic maps with values in a CAT(0) space under Lipschitz continuous boundary condition, due to Korevaar and Schoen in \cite{KorevaarSchoen1}, 
a result that extends the Hamilton theorem in \cite{Hamilton75}.

\br 
\label{remsimplyconnected}
Note that we cannot assume $Y$ to be only locally CAT(0). The fact that $Y$ is simply connected will be important here. 
Indeed, it is not clear 
how to parametrize the set of harmonic maps from the unit disk 
to a compact hyperbolic surface.  Similarly, it is not clear how to parametrize  the set of all harmonic maps from the unit disk to the circle $\m R/\m Z$, because this is equivalent to parametrizing 
all the harmonic functions on the unit disk.
\er

\br
Theorem \ref{thmbijectiveboundary} is an analog of the
theorems that parametrize unbounded harmonic maps between pinched Hadamard manifolds 
by their ``quasi-symmetric'' boundary condition at infinity. See the successive papers \cite{Markovic17}, \cite{LemmMarkovic},
\cite{BH15}, \cite{BH18}, and \cite{SidlerWenger}. 
that deal with an increasing level of generality.
\er

\subsection{Strategy of proof}
\label{secstrategy}

We will split the statement of  Theorem \ref{thmbijectiveboundary}
into five propositions.

\bp
\label{pronontangentiallimit}
{\bf (Construction of the boundary map) }
Let $ X$ be a \GGG Riemannian manifold
and let $Y$ be a proper {\rm CAT(0)}-space.
Let  $h: X\ra Y$ be a bounded harmonic map. 
Then, for $\si$-almost all $\xi\in \bX$,
the map $h$ admits a non-tangential limit $\ph(\xi):=\underset{x\ra\xi}{\rm NTlim}\,h(x)$ at the point $\xi$.

\ep

We denote by 
\begin{equation}
\label{eqnboundarymap}
\beta h:=\ph\in L^\infty(\bX,Y)
\end{equation}
the bounded measurable map from $\bX$ to $Y$
given by Proposition \ref{pronontangentiallimit}.
This map $\be h$ is called the boundary map of $h$, and the map 
$$
\beta:\mc H_b( X,Y)\ra L^\infty(\bX,Y)
$$
is called the boundary transform.

\bp
\label{proinjectiveboundary}
{\bf (Injectivity of the boundary transform) }
Same notation.
Two harmonic maps $h$, $h'$ from $ X$ to $ Y$
with  $\be h=\be h'$ are equal.
\ep

In order to prove that the transformation $\be$ is onto, we will construct its inverse map $P$.
We first rely on  the theorem
that solves the Dirichlet problem for harmonic maps with regular boundary data. It is due to Hamilton 
in \cite{Hamilton75} when the target is a manifold, and to Korevaar-Schoen in \cite[Thm 2.2]{KorevaarSchoen1} when the target is a  
{\rm CAT(0)} space.

\bfa 
\label{fachamilton}
{\bf (Hamilton, Korevaar and Schoen)}\\  Let $ \Om$ be a bounded Lipschitz Riemannian domain
and $Y$ be a complete {\rm CAT(0)} space.
Let $\ph:\partial \Om\ra Y$ be a Lipschitz map. 
Then, there exists a unique harmonic map $h=P\ph$ from $ \Om$ to $ Y$ that extends
continuously  $\ph$.
\efa

We then need to extend Fact \ref{fachamilton} to continuous boundary data, and to deal with a boundary   at infinity. This is included in the following
proposition wich will be proven in Section \ref{secdirichletproblem}. 
Note that when $X$ is a pinched Hadamard manifold and $Y$ is a 
CAT(0) Riemannian manifold, this  proposition is 
already in
\cite[Thm 3.2 and 4.7]{AvilesChoiMicallef91}.

\bp
\label{prodirichlet}{\bf (Dirichlet problem with continuous data)}
Let $ X$ be a \GGG Riemannian manifold 
and $Y$ be a complete {\rm CAT(0)}-space.
Let $\ph:\partial X\ra Y$ be a continuous map. 
Then, there exists a unique harmonic map $h=P\ph$ from $ X$ to $ Y$ that extends
continuously  $\ph$.
\ep

The conclusion in Proposition \ref{prodirichlet} means that the map $\ol h:\ol X\ra Y$ that is equal to $h$ on $ X$ and to $\ph$
on $\bX$ is continuous. Idem for Fact \ref{fachamilton}.

The main result in this article,
Theorem \ref{thmbijectiveboundary}, extends 
Proposition \ref{prodirichlet} to 
more general boundary conditions $\ph$. Indeed, it allows $\ph$ to be 
any bounded 
measurable map from $\partial X$ to $Y$.
As  will be very clear in the next proposition,
the proof of our main Theorem \ref{thmbijectiveboundary}
relies on the Hamilton, Korevaar, Schoen theorem.
\vs 

We endow  the space $L^\infty(\bX,Y)$ 
with ``the topology of the convergence in probability'', see \eqref{eqndphphp}.
The subspace $C(\bX,Y)$ of continuous maps  is then dense in  $L^\infty(\bX,Y)$, see Lemma \ref{lemcontinuousdense}.

We also endow the space $\mc H_b( X,Y)$ of bounded harmonic maps $h: X\ra Y$  with the topology of the uniform convergence on compact subsets of $ X$.

\bp
\label{propoissontransform}
{\bf (Construction of the Poisson transform) } Let $ X$ be a \GGG Riemannian manifold
and $Y$ be a bounded complete {\rm CAT(0)}-space.
The  map 
$$P:C(\bX,Y)\ra \mc H_b( X,Y)$$ given by Proposition $\ref{prodirichlet}$
has a unique continuous extension 
$$P:L^\infty(\bX,Y)\ra \mc H_b( X,Y).$$
\ep

We still call the extended map $P$ the Poisson transform. 

\bp
\label{prosurjectiveboundary}
{\bf (Surjectivity of the boundary transform) } Let $ X$ be a \GGG Riemannian manifold
and $Y$ be a compact {\rm CAT(0)}-space.
For all $\ph\in L^\infty(\bX,Y)$, one has $\ph=\be P\ph$.
\ep

\subsection{Overview}
\label{secplan}
\hspace{1em}

In Chapter \ref{secpreliminaryresult}, 
we recall preliminary facts about harmonic, subharmonic  and superharmonic functions $u$
on a \GGG Riemannian manifold. 
The key points that we will use are a control on the Poisson kernel in Proposition \ref{propoisson2},
upper bounds on the harmonic measures in Lemmas \ref{lemcontrolmeasure} and 
\ref{lemsidoubling},
and  the existence of non-tangential limit 
for bounded Lipschitz  superharmonic functions
on $ X$ in Proposition \ref{pronontangentiallimit2}. 

In Chapter \ref{secboundarymap}, we recall two facts about harmonic maps. 
The first one is the control, due to 
Cheng,  of the Lipschitz constant of a harmonic map
(Lemma \ref{lemcheng}).
The second one is  the subharmonicity of the distance function between two harmonic maps (Lemma \ref{lemharmonicsubharmonic}). We use these two facts, together with Proposition \ref{pronontangentiallimit2}, to prove 
the existence of non-tangential limit 
for our bounded harmonic map $h: X\ra Y$
(Proposition \ref{pronontangentiallimit}).
This provides the construction of the boundary map $\ph=\be h:\bX\ra Y$. These arguments also prove that the boundary transform $\be:h\mapsto \ph$ is injective (Proposition \ref{proinjectiveboundary}).  

In Chapter \ref{secpoissontransform}, we first 
construct the Poisson transform $P:\ph\mapsto h$
when the boundary data $\ph$ is continuous 
(Proposition \ref{prodirichlet}) by
building on the Hamilton, Korevaar 
and Schoen theorem (Fact \ref{fachamilton}). 
We then extend this transform $P:\ph\mapsto h$ to bounded measurable boundary data $\ph$ 
(Proposition \ref{propoissontransform}). The key point is a 
suitable uniform continuity property of this transform $\ph\mapsto h$. 

In Chapter \ref{secboundarypoisson}
we prove that the Poisson transform $P$ 
is a right inverse for the boundary transform $\be$, so that the boundary transform  $\be$ is surjective (Proposition \ref{prosurjectiveboundary}).
The key point is an estimate on 
sequen\-ces of subharmonic functions (Lemma \ref{lemradialconvergence}) that relies on  
the control of the Poisson kernel $P_\xi$ in 
Proposition \ref{propoisson2}
and on the Lebesgue density theorem for a doubling measure on a compact quasi-metric space
(Fact \ref{faclebesgue}).
\vs 

We would like to thank F. J\"{a}ckel for useful comments on a first version of this preprint.

\section{Harmonic and subharmonic functions}
\label{secpreliminaryresult}
\bq 
In this second chapter, we gather a few results concerning harmonic and subharmonic functions on a \GGG Riemannian manifold $X$ that will be used in the proof ot Theorem \ref{thmbijectiveboundary}. 
\eq 
We set $g$ for the Riemannian metric,
$d$ for the Riemannian distance,
$\Delta$ for the Laplace Beltrami operator  and $k=\dim X$. For the Potential theory of a \GGG Riemannian manifold, we refer 
to the seminal paper \cite{AnconaStFlour} 
and to its recent update \cite{KemperLohkamp}.

\subsection{The Harnack inequality and the Green function}
\label{secHarnack}
\bq
In this section we present three  classical Harnack inequalities for positive harmonic functions.
\eq

We recall that a function $u: X\ra \m R$ is superharmonic
if it  is lower semicontinuous, locally integrable, and if $\Delta u\leq 0$ holds in the weak sense. 
A function $u$ is subharmonic if $-u$ is superharmonic. 
A function $u$ is  harmonic if 
it is both subharmonic and superharmonic

The Harnack inequality, which has been improved by Serrin and by S.T. Yau, gives a uniform control for positive harmonic functions
on compact sets.
See \cite[Lemma 2.1]{LiWang02} for a short proof, and also \cite[Cor. 8.21]{GilbargTrudinger}.

\bfa
\label{facHarnack1}{\bf (Harnack inequality)} 
Let $X$ be a complete Riemannian manifold with bounded sectional curvature.
There exists a constant $c_0>0$  such that, for any positive harmonic function $u$ on a
ball 
$
B(x_0, r)
$
with $r\leq 1$,
one has 
$$
\|\nabla \log u(x)\|\;\leq\; c_0/r
\;\;\;\mbox{\rm for all $x$ in $B(x_0,r/2)$}.
$$
One then has 
$$
u(y)\;\leq\; e ^{c_0}u(x)
\;\;\;\mbox{\rm for all $x$, $y$ in $B(x_0,r/2)$}.
$$
\efa

The Green operator $G$ is the ``inverse'' of 
the Laplacian. The spectral gap assumption ensures that 
the Green operator is bounded as an operator on $L^2(X)$.
The Green kernel $G(x,y)$ is the kernel of the Green
operator. It is is symmetric i.e. 
$G(x,y)=G(y,x)$. It is a positive $C^\infty$-function
on $X\times X\smallsetminus \De_X$ and, for each $x$
in $X$, the function $G_x:=G(x,.)$ satisfies 
$\De G_x=-\de_x$. In particular the function $G_x$ is harmonic 
outside $\{x\}$.

The bounded geometry assumption ensures the following control on the Green function. We set $\log_*(t):={\rm max}(1, \log t)$.

\bfa
\label{facGreen}{\bf (Green function)} 
Let $X$ be a Riemannian manifold with bounded geometry and with spectral gap. There exist $C_0>1$ and $\eps_0>0$ 
such that~:\\
$a)$   For $x, y$ in $ X$ with 
$d(x,y)\leq 1$, one has
\begin{eqnarray*} 
\label{eqngreenproche}
C_0^{-1}d(x,y)^{2-k}\leq G(x,y)\leq C_0\,d(x,y)^{2-k}
\hspace{1.2em}
& &
\mbox{if $k\neq 2$,}\\
\nonumber
C_0^{-1}\log_*(1/d(x,y))\leq G(x,y)\leq C_0\,\log_*(1/d(x,y))
\hspace{-1em}
& &
\mbox{if $k= 2$.}
\end{eqnarray*}
$b)$ For $x, y$ in $ X$ with 
$d(x,y)\geq 1$, one has
\begin{eqnarray} 
\label{eqngreenloin}
G(x,y)&\leq &C_0 \, e^{-\eps_0\, d(x,y)}
\end{eqnarray}
\efa

See for instance \cite[Prop. 2.7 and 2.12]{KemperLohkamp}.

\subsection{The Ancona Inequality}
\label{secAncona}

The Gromov hyperbolicity assumption ensures a much more precise control on the Green function due to Ancona.

\bfa
\label{facAncona}
{\bf (Ancona Inequality)} 
Let $X$ be a \GGG Riemannian manifold. 
Then, there exists $C_1>1$ such that for any  point $y$
on a geodesic segment $[x,z]$ in $X$  
such that $d(x,y)\geq 1$ and $d(y,z)\geq 1$
one has
\begin{equation}
\label{eqnAncona}
C_1^{-1}\,G(x,y)\,G(y,z)\;\leq\; G(x,z)\; \leq\;
C_1\,G(x,y)\,G(y,z)
\end{equation}
\efa

The boundary Harnack inequality
compares the behavior of two positive harmonic functions
near a piece of the boundary $\bX$ 
where they both go to zero.
In order to state this inequality we need to 
introduce some notation. 

We first recall the definition of the Gromov product 
for two  points $\eta_1$, $\eta_2$ in $\ol X=X\cup \bX$
seen from a point $o\in X$~:
$$
(\eta_1|\eta_2)_{o}:=
\limsup\limits_{\substack{x_1\ra\eta_1\\x_2\ra\eta_2}}
(x_1|x_2)_{o}.
$$
This quantity is equal, up to a uniformly bounded 
error term, to the distance between 
$o$ and a geodesic going from $\eta_1$ to $\eta_2$.

For $x$ in $\ol X$, we introduce the sets
\begin{equation*}
\mc H_o^{t}(x)=\{y\in  X\mid (y|x)_o\geq t\}, 
\end{equation*}
\begin{equation*}
\ol{\mc H}_o^{t}(x)=\{y\in \ol X\mid (y|x)_o\geq t\}.
\end{equation*}
Note that these sets are empty when $d(o,x)<t$.

We recall that the topology of $\ol X$ is the topology that extends the topology of $X$ 
and such that a neighborhood basis of a point $\xi\in \bX$
is given by the sets $\ol{\mc H}_o^{t}(\xi)$ with $t>0$.
See \cite{GhysHarp90}.
We will call them (rough) half-spaces.

We also recall that, since $X$ is Gromov hyperbolic, 
we can choose  $\de>1$ satisfying  \eqref{eqngromovproduct} and such that, for every 
geodesic triangle  $x$, $y$, $z$ in $X$, 
every point $u$ on the edge $xz$ 
is at distance at most $\de$ of the union  $xy\cup yz$ of the other 
edges.

We record a five properties of these half-spaces~:
\begin{equation}
\label{eqnhsxyhtyx}
\ol{\mc H}_x^{s}(y)\cup \ol{\mc H}_y^{t}(x)
\; =\; \ol X
\;\; 
\mbox{\rm for all $x$, $y$ in $X$  with 
$d(x,y)\geq s+t$,}
\end{equation}
\begin{equation}
\label{eqnhsxyhtyxbis}
\ol{\mc H}_x^{s}(y)\cap \ol{\mc H}_y^{t}(x)
\; =\; \emptyset
\;\; 
\mbox{\rm for all $x$, $y$ in $X$  with 
$d(x,y)< s+t$.}
\end{equation}
\begin{equation}
\label{eqnhtxyhtxz}
\ol{\mc H}_x^{t}(y)\subset \ol{\mc H}_x^{t-\de}(z)
\;\; 
\mbox{\rm for all $x$ in $X$, $y$, $z$ in $\ol X$  with 
$z\in\ol{\mc H}_x^{t}(y)$,}
\end{equation}
\begin{equation}
\label{eqnhtoyhtmy1}
\ol{\mc H}_x^{t+s}(y)\subset \ol{\mc H}_m^{t}(y)
\;\; 
\mbox{\rm for all $m$, $x$, $y$ in $X$ with $d(x,m)=s$ ,}
\end{equation}
\begin{equation}
\label{eqnhtoyhtmy}
\ol{\mc H}_m^{t}(y)
\subset \ol{\mc H}_x^{t+s-\de}(y)
\;\; 
\mbox{\rm for $m$ in $xy$ with $d(x,m)=s$ and $t>\de$.} 
\end{equation}
These properties explain why these sets are called 
half-spaces. 

The first four ones are straightforward. 
For the last one, notice that, since $t>\de$, 
for any $z$ in $\ol{\mc H}_m^{t}(y)$, the distance between $m$ 
and a geodesic $yz$ is larger than $\de$ and hence 
there exists a point $m'$ on a geodesic $xz$ such that
$d(m,m')\leq \de$.
\vs

We now state the strong boundary Harnack inequality for
a  \GGG Riemannian manifold $X$. This inequality 
is actually equivalent to the Ancona Inequality 
in Fact \ref{facAncona}.

\bfa
\label{facHarnack3}{\bf (Strong Boundary Harnack inequality)} 
Let $X$ be a \GGG Riemannian manifold.
There exists $C_2>0$ and $t_0>0$ such that 
for all $o\in X$, all $\xi\in \bX$, all $t\geq 0$ 
and  all positive continuous functions $u$, $v$  on the 
half-space $\ol{\mc H}:=\ol{\mc H}^{t}_{o}(\xi)$
which are zero on $\bX\cap \ol{\mc H}$ and harmonic in the interior of  $X\cap \ol{\mc H}$,
one has
\begin{eqnarray} 
	\label{eqnHarnack3}
	\frac{u(y)}{v(y)}\leq C_2\, \frac{u(x)}{v(x)}
	&\mbox{\it for all}& x, y \in \mc H^{t+t_0}_o(\xi).
\end{eqnarray}
\efa

\br 
Fact \ref{facHarnack3} is due to Anderson and Schoen in \cite[Corollary 5.2]{AndersonSchoen} when $X$ is a pinched Hadamard manifold. It is due to Ancona 
in \cite{AnconaStFlour} when $X$ is a 
\GGG Riemannian manifold.
\er

The following statement is a Corollary 
of Facts \ref{facHarnack1} and \ref{facHarnack3}.

\bc
\label{corHarnack4}{\bf (Boundary Harnack inequality)} 
Let $X$ be a \GGG Riemannian manifold.
Let $O\subset \ol X$ be a connected open subset and $K\subset O$ be a compact subset. There exists a constant $C=C_{K,O,X}>0$
such that for all continuous functions 
$u$, $v$ on $O$ that are 
harmonic and positive on $O\cap  X$ and zero on $O\cap \bX$, 
one has
\begin{eqnarray} 
	\label{eqnHarnack4}
	\frac{u(y)}{v(y)}\leq C\, \frac{u(x)}{v(x)}
	&\mbox{\it for all}& x, y \in K\cap X.
\end{eqnarray}
\ec

Note that   Fact \ref{facHarnack3} is stronger 
than Corollary \ref{corHarnack4} since it requires the constant $C_2$ 
not to depend on the half-space $\ol{\mc H}$.\vs

\subsection{The Poisson kernel}
\label{secpoissonfunction}

\bq
We now recall the definition of the Poisson kernel,
also called the Martin kernel. 
\eq

The following fact is due to 
Anderson and Schoen 
in \cite{AndersonSchoen} for a pinched Hadamard manifold,
and has been generalized by Ancona in \cite{AnconaStFlour}
for a \GGG Riemannian manifold $X$. See also \cite{KemperLohkamp} for more general measured metric spaces.
It describes all the positive harmonic functions on $X$ 
and, more precisely, it describes the Martin boundary of $X$.
The key point in the proof is the strong boundary Harnack inequality of Fact \ref{facHarnack3}.

\bfa {\bf (Martin Boundary)}
\label{facMartin}
Let $X$ be a \GGG Riemannian manifold, and fix a point $o\in X$.\\
$a)$ For all $\xi$ in $\bX$, 
there exists a unique non-negative continuous function
$$
x\mapsto P_\xi(x)=P_\xi(o,x)
$$ 
on $\ol X\smallsetminus\{\xi\}$ 
which is harmonic on $X$, zero on $\bX \smallsetminus\{\xi\}$ with $P_\xi(o)=1$.\\ 
$b)$ For $x\in X$ and $\xi\in \bX$, $P_\xi(x)$ is obtained as the limit
\begin{equation}
P_\xi(x)=\lim_{y\ra\xi}\,\frac{G(x,y)}{G(o,y)}\, .
\end{equation}
$c)$ Any positive harmonic function $h$ on $X$ can be written as
$$
h(x)=\int_{\bX}P_\xi(x)\,{\rm d}\mu(\xi)
$$ 
for a unique 
positive finite Borel measure $\mu=\mu_h$ on $\bX$.\\
$d)$ The function $(x,\xi)\mapsto P_\xi(x)$ 
is continuous on $\ol{X}\times\bX\smallsetminus \De_{\bX}$ where $\De_{\bX}$ denotes the diagonal in $\bX\times\bX$.
\efa
By the very definition of the Poisson functions, the following holds~:
\begin{equation*}
\label{eqnpoissonsymmetric}
P_\xi(x,o)=P_\xi(o,x)^{-1}.
\end{equation*}
Those functions $P_\xi$ are exactly the positive harmonic functions on $X$
that are minimal, up to normalization: i.e. the only positive harmonic functions $h$ on $X$ that are bounded by $P_\xi$ are the multiples $h=\al P_\xi$ with $\al\leq 1$.

In the proof of Proposition \ref{prosurjectiveboundary}, 
we will need the following estimate for the Poisson functions.

\bp
\label{propoisson2} 
Let $X$ be a \GGG Riemannian manifold.
Then, there exists a constant $C_3>0$ such that for all
$o$ in $X$, all $\xi$ in $\bX$, all $x$ on a geodesic ray $o\xi$, and 
for all $\eta_1$, $\eta_2$ in $\bX$
with 
\begin{equation}
\label{eqnhyppoisson2}
|(\xi|\eta_1)_{o}-(\xi|\eta_2)_{o}|\leq 1\, ,
\end{equation}
one has
\begin{equation}
\label{eqnpoisson2}
C_3^{-1}\leq \frac{P_{\eta_{_1}}(x)}{P_{\eta_{_2}}(x)}\leq C_3\, .
\end{equation}
\ep

We begin by a lemma that extends  Ancona's inequality to 
Poisson functions.
Since $X$ is Gromov hyperbolic, for every 
geodesic triangle with distinct vertices $x$, $y$, $z$ in $\ol X$, there exists a point $m$ 
whose distance to any of the three geodesic sides is at most $\de$.
Such a point $m$ is called a center of the triangle $x,y,z$.
It is not unique, but the distance between two centers is 
at most $8\de$. See \cite{GhysHarp90}.

\bl
\label{lemAncona2}
Let $X$ be a \GGG Riemannian manifold. 
Then, there exists a constant  $C_4>1$ such that 
for all $\xi$ in $\bX$, all points $o$, $x$ in $X$ 
one has
\begin{equation}
\label{eqnAncona2}
C_4^{-1}\,\frac{G(m,x)}{G(m,o)}\;\leq\; P_\xi(o,x)\; \leq\;
C_4\,\frac{G(m,x)}{G(m,o)},
\end{equation}
where $m$ is a center of a geodesic triangle with vertices
$o$, $x$, $\xi$ such that $d(o,m)\geq 1$ and $d(x,m)\geq 1$.
\el

Note that it is possible to choose such a center $m$ since $X$ is Gromov  hyperbolic with constant $\de>1$.
\vs 

A particular instance of \eqref{eqnAncona2},  when $x$ is on a geodesic ray from $o$ to $\xi$ and $d(o,x)\geq 1$, reads as  
\begin{equation}
\label{eqnAncona2bis}
C_4^{-1}G(o,x)^{-1}\;\leq\; P_\xi(o,x)\; \leq\;
C_4\, G(o,x)^{-1}. 
\end{equation}

\begin{proof}
The center $m$ is at a distance at most $\de$ from both a geodesic
ray $o\xi$ and a geodesic ray $x\xi$. 
Therefore when a point $y\in X$ is sufficiently near $\xi$,
the point $m$ is  at distance at most 
$2\de$ from both the geodesic rays $oy$ and $xy$.
Therefore, using the Harnack inequality in Fact 
\ref{facHarnack1} and the 
Ancona inequality in Fact \ref{facAncona} one gets,
with a constant $c_4$ depending only on $X$~:
\begin{eqnarray*}
\label{eqnAncona3}
c_4^{-1}\,G(x,m)G(m,y)\;\leq\; G(x,y)\; \leq\;
c_4\,G(x,m)G(m,y),\\
c_4^{-1}\,G(o,m)G(m,y)\;\leq\; G(o,y)\; \leq\;
c_4\,G(o,m)G(m,y).
\end{eqnarray*}
Taking the ratio of these estimates yields, with $C_4=c_4^2$,
\begin{equation*}
\label{eqnAncona4}
C_4^{-1}\,\frac{G(x,m)}{G(o,m)}\;\leq\; \frac{G(x,y)}{G(o,y)}\; \leq\;
C_4\,\frac{G(x,m)}{G(o,m)}.
\end{equation*}
One then gets \eqref{eqnAncona2} by letting $y$ converge to $\xi$.
\end{proof}

\begin{proof}[Proof of Proposition \ref{propoisson2}] 
Let $m_1$ be a center of the triangle $o$, $x$, $\eta_1$
and $m_2$ be a center of the triangle $o$, $x$, $\eta_2$.
Those points can be chosen so that $d(m_i,o)\geq 1$ 
and $d(m_i,x)\geq 1$.
The assumption \eqref{eqnhyppoisson2} tells us that the distance $d(m_1,m_2)$ is bounded by a constant depending only on $\de$.
Therefore by the Harnack principle in Fact \ref{facHarnack1},
there exists a constant $c_3>1$ depending only on $X$ such that
\begin{equation}
\label{eqnAncona5}
c_3^{-1}\leq \frac{G(x,m_2)}{G(x,m_1)}\leq c_3
\;\;{\rm and}\;\; 
c_3^{-1}\leq \frac{G(o,m_2)}{G(o,m_1)}\leq c_3
\end{equation}
Taking the ratio of the bound \eqref{eqnAncona2} with $\xi=\eta_1$ 
by the bound \eqref{eqnAncona2} with $\xi=\eta_2$ 
one gets
\begin{equation}
\label{eqnpoisson3}
C_4^{-2}\,\frac{G(x,m_1)}{G(x,m_2)}\,\frac{G(o,m_2)}{G(o,m_1)}
\;\leq\; 
\frac{P_{\eta_{_1}}(x)}{P_{\eta_{_2}}(x)}
\;\leq\; 
C_4^2\,\frac{G(x,m_1)}{G(x,m_2)}\,\frac{G(o,m_2)}{G(o,m_1)}.
\end{equation}
Combining  \eqref{eqnAncona5} with \eqref{eqnpoisson3}, one obtains 
\eqref{eqnpoisson2} with $C_3=c_3^2C_4^2$.
\end{proof}

\bc
\label{corAncona2}
Let $X$ be a \GGG Riemannian manifold. 
Then, there exists a constant  $C_4>1$ such that 
for all $\xi$ in $\bX$, all points  $x$ in $X$ 
and all point $y$ on a ray $x\xi$ with $d(x,y)\geq 1$,
one has
\begin{equation}
\label{eqnAncona6}
{C_4}^{-1}\, G(x,y)\,P_\xi(y)\;\leq\;
P_\xi(x)
\; \leq\;
C_4\,G(x,y)\,P_\xi(y).
\end{equation}
\ec
\begin{proof}
Since, by definition $P_\xi(x,y)=P_\xi(y)/P_\xi(x)$,
inequalities \eqref{eqnAncona6} are nothing but a reformulation of \eqref{eqnAncona2bis}.
\end{proof}

See also \cite[Cor. 6.4]{AndersonSchoen} and \cite{LedrappierLim}  for other estimates on the 
Poisson kernel when $X$ is a pinched Hadamard manifold.

We recall that, in Geometric Group Theory, the word ``geodesic'' means 
``minimizing geodesic''.  This is why  the following lemma is non-trivial.

\bl
\label{lemnodeadend}
Let $X$ be a \GGG Riemannian manifold. 
Then there exists  $C_5>1$ such that, 
for all $x$ in $X$, there exists $\xi$, $\eta$ in $\bX$
with $(\xi|\eta)_x\leq C_5$.
\el

We do not explicitely use this lemma. 
It illustrates the influence of the spectral gap condition
on the geometry of $X$. It tells us there are no dead ends in $X$.
  Here is a sketch of proof.

\begin{proof}
Since $X$ is Gromov hyperbolic,
if there were dead ends, 
we would be able to find, for all $n\geq 1$, 
a ball $B_n:=B(x_n, n)$ of radius $n$ in $X$
whose boundary  $S_n:=S(x_n,n)$ has diameter at most $2\de$.
Since $X$ has bounded geometry, the volumes of the balls
$B_n$ go to infinity while the diameters
of $S_n$ are bounded.
This contradicts the spectral gap.
\end{proof}

\subsection{The harmonic measures}
\label{secharmonicmeasure}

\subsubsection{Harmonic measures for a bounded domain}

We first recall the solution of the Dirichlet problem 
for harmonic functions on a Lipschitz bounded Riemannian domain 
$\Om$ (see Section \ref{seclipschitzdomain} for a precise definition). This can be found in \cite[Chapter 6 and 8]{GilbargTrudinger} 
when $\bOm$ is smooth and in \cite{Ancona78} when $\bOm$ is Lipschitz continuous. It says:

\bfa {\bf (Dirichlet problem for functions on a bounded domain)} Let $\Om$ be a Lipschitz bounded Riemannian domain.
 For every continuous function $\ph\in C(\bOm,\m R)$, there exists a unique continuous function $h:\ol\Om\ra \m R$ which is harmonic on $\Om$ and equal to $\ph$
on $\bOm$. 

Moreover,
for $x\in  \Om$  there exists a measure $\si_x=\si^\Om_x$  on $\bOm$, called the harmonic measure on $\bOm$ seen from $x$, 
such that the harmonic extension $h$ of $\ph$ is given by 
\begin{eqnarray}
\label{eqnuintphiom}
h(x):=\int_{\bOm}  \ph(\xi)\,{\rm d}\si_x(\xi).
\end{eqnarray}
\efa
\noindent
By a theorem of Dahlberg in \cite{Dahlberg79},
the harmonic measure $\si^{\Om}_{x}$ at each point $x \in\Om$ is equivalent to the 
Riemannian measure  on $\bOm$. 

The measure  $\si_x^\Om$ is  a doubling measure on $\bOm$.
This means, see \cite[Section 11.3]{CaffarelliSalsa},  that there exists a constant 
$c=c_{\Om,x}$ such that for all $r>0$ and 
all $\xi$ in $\bOm$ one has 
\begin{equation}
\label{eqndoubling}
\si_{x}^\Om(B(\xi,2r))\leq c\, \si_{x}^\Om(B(\xi,r)).
\end{equation}
In this notation, we think of $\si_x^\Om$ as a measure on $X$
supported by $\bOm$.
\br
From a probabilistic point of view, the harmonic measure $\si^\Om_x$ on $\bOm$
is the exit probability measure
of a Brownian motion on $ \Om$ starting at point $x$.
\er

\subsubsection{Harmonic measures on \GGG Riemannian manifolds }
We now recall the solution  of the Dirichlet problem 
for harmonic functions on a \GGG Riemannian manifold $X$. 
This is independently due to Anderson and Sullivan
when $X$ is a pinched Hadamard manifold.
See \cite{AndersonSchoen} for a nice account. 
It is due to Ancona when $X$ is a \GGG Riemannian manifold, 
as a consequence of the description of the Martin boundary of $X$. 

\bfa
{\bf (Dirichlet problem for functions on \GGG manifolds)}\\
Let $X$ be a \GGG Riemannian manifold.
For every continuous function $\ph\in C(\bX,\m R)$, there exists a unique continuous function $h:\ol X\ra \m R$ which is harmonic on $ X$ and is equal to $\ph$
on $\bX$. 

Moreover,
for $x\in  X$  there exists a measure $\si_x=\si^X_x$  on $\bX$, called the harmonic measure on $\bX$ seen from $x$, 
such that the harmonic extension $h$ of $\ph$ is given by 
\begin{eqnarray}
\label{eqnuintphix}
h(x):=\int_{\bX}  \ph(\xi)\,{\rm d}\si_x(\xi).
\end{eqnarray}
\efa 

For $x=o$, the probability measure $\si_o$ on $\bX$ is the one
that appears in 
the decomposition of the constant harmonic function $h=1$ in Fact \ref{facMartin}.$c$. 

For every $x$ in $X$, the positive measure $\si^X_x$ 
is given by the formula
\begin{equation*}
{{\rm d}\si^X_x(\xi)=	P_\xi(x)\,{\rm d}\si^X_{o}}(\xi)
\end{equation*}
so that Equation \eqref{eqnuintphix}  can be rewritten as
\begin{equation}
\label{eqnuintphi2}
h(x):=\int_{\bX}  \ph(\xi)\,P_\xi(x)\,{\rm d}\si^X_{o}(\xi).
\end{equation}

When $\ph$ is continuous,  the function $h$ 
defined on $X$ by \eqref{eqnuintphi2} is harmonic and extends continuously $\ph$. 
Indeed each function $\xi\mapsto P_\xi(x)$, for ${x\in  X}$, 
is positive and satisfies
$\int_{\bX}P_\xi(x)\,{\rm d}\si_{o}^X(\xi) =1$ and, when a sequence $(x_n)$ converges to $\xi\in \bX$, the sequence of probability measures $(P_\xi(x_n){\rm d}\si^X_{o}(\xi))$  converges weakly to $\de_{\xi}$.

Note that, even when $X$ is a pinched Hadamard manifold,
the measure $\si_{o}^X$ is not always equivalent 
to the ``visual measure''.
\vs 

In order to have shorter notation, 
we will think of the harmonic measures $\si_x^X$ as 
measures on $\ol X$ supported by $\bX$.

\subsubsection{Upper bound for the harmonic measures}
\label{secupperhar}

We will need the following  
uniform control of the harmonic measures
$\si^{\Om}_{o}$
for bounded subdomains of $X$ with Lipschitz boundary. By definition, the
probability measure $\si^{\Om}_{o}$ is supported by the boundary $\bOm$.
This control tells us that, seen from $o$,
the measure of the part of  $\bOm$
cut out by a half space far away from $o$ is uniformly small.

\bl 
\label{lemcontrolmeasure} 
Let $X$ be a \GGG Riemannian manifold.
For all $\eps >0$ there  exists $\ell=\ell_\eps>0$ such that
for   all $o$ in $X$, $x$ in $\ol X$, 
one has 
\begin{equation}
\label{eqncontrolx}
\si^{X}_{o}(\ol{\mc H}^{\ell}_{o}(x))
\;\leq\; \eps,
\end{equation}
and, for all bounded Lipschitz subdomain $\Om\subset X$ containing $o$, one has 
\begin{equation}
\label{eqncontrolom}
\si^{\Om}_{o}(\mc H^{\ell}_{o}(x))
\;\leq\; \eps.
\end{equation}
\el

\begin{proof} We first prove \eqref{eqncontrolom}.
We introduce the set 
$$
E:=\mc H^{\ell}_{o}(x)\cap \bOm ,
$$
where $\ell>\de$ will be chosen later, and the open $1$-neighborhood  of $E$
$$U:=\{y\in X\mid d(y,E)<1\}.
$$

We introduce then the reduced function $u:=R_1^{U}$ of the constant function $1$ to this open set $U$. 
By definition, $u$ is the smallest positive superharmonic function on $X$ which is larger than ${\bf 1}_U$. This function is equal to $1$ on $U$, 
it is harmonic on $X\smallsetminus\ol U$
and one has $0\leq u\leq 1$. Since $U$ is relatively compact 
this function $u$ is a {\it potential} on $X$ i.e. its largest harmonic minorant is $0$.

The Riesz decomposition theorem tells us that every potential $u$ on $X$ can be written in a unique way as 
\begin{equation}
\label{eqnbalayage}
u(x)=\int_X G(x,y)\,{\rm d}\la(y)\, ,
\end{equation}
where $\la$ is a positive Radon measure
on $X$ 
called the {\it Riesz measure} of $u$.

In our case where $u$ is the reduced function $u=R^U_1$ for a relatively compact 
open set $U$,
the Riesz measure $\la$ is a finite measure supported by the boundary $\partial U$.  
\vs 
 
Since $u$ is a positive superharmonic function on $X$ 
which is equal to $1$ on $E$
one has, for all $z\in \Om$,
\begin{equation}
\label{eqnsibxux}
\si^{\Om}_{z}(E)\leq u(z).
\end{equation}

We can assume that $E$ is not empty, and hence that 
$d(o,x)\geq \ell$. Let $m$ be a point on a geodesic segment from $o$ to $x$ 
with $d(o,m)=\ell$.
We claim that there exists a constant $C>0$
depending only on $X$  such that
\begin{eqnarray}
\label{eqnuxuxgx}
u(o)\leq C\, G(o,m) \, u(m).
\end{eqnarray}

Indeed, 
for each $y$ in $\mc H^{\ell}_{o}(x)$, any  geodesic segment from $o$ to $y$ intersects the 
ball $B(m,\de)$.
Applying Harnack Inequality and Ancona Inequality, 
one finds a constant $C_6>0$ depending only on $X$ such that, for all $t>1$ and all $y$ in $\partial U$, one has 
$$
G(o,y)\leq C_6\, G(o,m)\, G(m,y).
$$

Applying this inequality to each of the Green functions
in the integral \eqref{eqnbalayage},
one gets our claim \eqref{eqnuxuxgx}.

Since $u$ is bounded by $1$, it follows from
\eqref{eqnsibxux}, \eqref{eqnuxuxgx} and  \eqref{eqngreenloin} that
$$
\si^{\Om}_{o}(E)
\leq u(o)\leq C_6\, G(o,m)
\leq C_6\, C_0\, e^{-\eps_0 \ell}\leq 
\eps
$$
if $\ell =\ell_\eps$ is chosen large enough.
\vs 

We now prove \eqref{eqncontrolx} with the same $\eps$ 
and $\ell_\eps$ as in  \eqref{eqncontrolom}.
If \eqref{eqncontrolx} were not true, there would exist a  point $x\in\ol X$, a small constant $\al>0$ and a continuous function $\ph:\bX\ra [0,1-\al]$ supported by an open subset of $\bX$ included in 
$\ol{\mc H}^{\ell}_{o}(x)$ whose harmonic extension
$h:X\ra [0,1-\al]$ satisfies $h(o)>\eps$.
Since $h$ is continuous, if $\Om$ contains a sufficiently large ball $B(o,R)$, the restriction of $h$ to
the complement
$\bOm \smallsetminus\mc H^{\ell}_{o}(x)$ 
is bounded by $\al\eps$.
Therefore, applying formula \eqref{eqnuintphiom} with $\ph=h$ and
using \eqref{eqncontrolom}, we get
$$
\eps
\;<\; 
h(0)
\;\leq\; 
(1-\al)\,\si_o^\Om({\mc H}^{\ell}_{o}(x))+\al\eps
\;\leq\; 
(1-\al)\,\eps +\al \eps 
\;=\;
\eps.
$$
This contradiction proves that \eqref{eqncontrolx} is true.
\end{proof}

\subsubsection{Lower bound for the harmonic measure}
\label{seclowerhar}

In order to prove the doubling property of the harmonic measure on $X$, we will need the following uniform lower bound 
on the harmonic measure

\bl 
\label{lemcontrold}
Let $X$ be a \GGG Riemannian manifold. 
For all $\ell\geq 0$, there exists $\eps_\ell>0$ 
such that
such that for all $o$ in $X$ and $\xi\in \bX$, one has
\begin{equation}
\label{eqncontrold}
\si_{o}^X(\ol{\mc H}^{\ell}_o(\xi))\;\geq\; \eps_\ell.
\end{equation}
\el

\begin{proof} Let $\ell_0>0$ be the length given by \eqref{eqncontrolx} 
with $\eps =1/2$.
Let $m$ be a point on a geodesic ray $o\xi$ such that 
$d(o,m)=\ell+\ell_0+\de$, so that, by  
\eqref{eqnhtxyhtxz} and \eqref{eqnhsxyhtyx},
$$
\ol{\mc H}^{\ell}_o(\xi)\supset 
\ol{\mc H}^{\ell+\de}_o(m)\supset
\ol X\smallsetminus \ol{\mc H}^{\ell_{_0}}_m(o).
$$ 
By the Harnack inequality applied to the harmonic function 
$x\mapsto \si_x^X(\ol{\mc H}^\ell_o(\xi))$, there exists
a constant $C_\ell>0$ depending only on $X$ and $\ell$ such that
\begin{equation*}
\si_o^X(\ol{\mc H}^{\ell}_o(\xi))
\;\geq\; 
C_\ell^{-1}\si_m^X(\ol{\mc H}^{\ell}_o(\xi))
\;\geq \; C_\ell^{-1}\,(1-\si_m^X(\ol{\mc H}^{\ell_{_0}}_m(o))).
\end{equation*}
The choice of $\ell_0$ implies
$\si_m^X(\ol{\mc H}^{\ell_{_0}}_m(o))\leq 1/2$.
This gives \eqref{eqncontrold} with the constant $\eps_\ell=1/(2C_\ell)$.
\end{proof}

\subsubsection{Doubling for the harmonic measure}
\label{secdoublinghar}

The following Lemma \ref{lemsidoubling} tells us that the measure $\si_o^{X}$ satisfies a doubling property.
See  \cite[Lemma 7.4]{AndersonSchoen} 
when $X$ is a pinched Hadamard manifold.

\bl 
\label{lemsidoubling}
Let $X$ be a \GGG Riemannian manifold. 
There exists a constant 
$c=c_{X}$ 
such that for all $o$ in $X$, $\xi\in \bX$ and $t\geq 0$, one has
\begin{equation}
\label{eqndoubling2}
\si_{o}^X(\ol{\mc H}^{t}_o(\xi))\;\leq\; 
c\, \si_{o}^X(\ol{\mc H}^{t+1}_o(\xi)).
\end{equation}
\el 

\begin{proof} By Lemma \ref{lemcontrold}, we may assume that 
$t\geq 1$.

We first claim that there exists a constant  $C_7>1$ such that for all $o$ in $X$, all $\xi$ in $\bX$, and all  $x_t$ on a geodesic ray $o\xi$ with $d(o,x_t)=t\geq 1$, one has
\begin{equation}
\label{eqnsiohtoxi} 
C_7^{-1}\, G(o,x_t)
\;\leq\;
\si_{o}^X(\ol{\mc H}^{t}_o(\xi))
\;\leq\;  
C_7 \, G(o,x_t)\, .
\end{equation}
In order to prove this claim, we  introduce for each $t>0$ the harmonic function 
$$
z\mapsto h_t(z):=\si^X_z(\ol{\mc H}^{t}_o(\xi))=
\int_{\ol{\mc H}^{t}_o(\xi)\cap\bX}P_\eta(z)\,{\rm d}\si_o^X(\eta).
$$
Integrating 
the inequalities \eqref{eqnAncona6}, one finds a constant $C_4>1$ depending only on $X$
such that 
\begin{equation}
\label{eqnhxthoxht}
{C_4}^{-1}G(o,x_t)h_t(x_t)
\;\leq\;
h_t(o)
\;\leq\;
{C_4}\,G(o,x_t)h_t(x_t).
\end{equation}
We recall from \eqref{eqnhtoyhtmy}
that 
\begin{equation*}
\label{eqnhotxi}
\ol{\mc H}^{2\de}_{x_{_t}}(\xi)
\;\subset\;
\ol{\mc H}^{t}_o(\xi),
\end{equation*}
so that, using Lemma \ref{lemcontrold} with $\ell=2\de$, 
one gets a constant $\eps_{2\de}>0$ such that
\begin{equation}
\label{eqnsixthxtd}
\eps_{2\de}
\;\leq\;
\si_{x_{_t}}(\ol{\mc H}^{2\de}_{x_{_t}}(\xi))
\;\leq\;
h_t(x_{t})
\;\leq\;
1.
\end{equation}
Combining \eqref{eqnhxthoxht} with \eqref{eqnsixthxtd},
we obtain our claim \eqref{eqnsiohtoxi}.

Now, the 
Harnack inequality in Fact \ref{facHarnack1}
provides a constant $C$ 
depending only on $X$ such that, for $t\geq 1$,
\begin{equation}
\label{eqngoxgoxgox}
G(o,x_{t})
\;\leq\;
C\, G(o,x_{t+1}). 
\end{equation}
The bound \eqref{eqndoubling2} follows from 
\eqref{eqnsiohtoxi} and \eqref{eqngoxgoxgox}.
\end{proof}

\subsection{Non-tangential limits}
\label{secradiallimit}

According to Fatou's theorem, every bounded harmonic function on 
the Euclidean ball $B\subset \m R^k$ admits 
a non-tangential limit at $\si_0$-almost all points 
of the boundary sphere $\bB$ (see \cite[Theorem 4.6.7]{ArmitageGardinerPotential}).
This is not always  true for a bounded superharmonic function $u$, see \cite[p. 175]{Tsuji75}.

Yet, according to Littlewood's theorem,
every bounded superharmonic function $u$ on 
$B$ admits a radial limit at $\si_0$-almost all points of $\bB$, 
see  \cite[Thm. 4.6.4, Cor. 4.6.8]{ArmitageGardinerPotential}.
One needs an extra assumption on $u$ to ensure that
this radial limit is also a non-tangential limit.
This  condition is the ``Lipschitz continuity of $u$
for the hyperbolic metric on the ball $B$''.

\bp
\label{pronontangentiallimit2} Let $X$ be 
a \GGG Riemannian manifold, $o\in X$ and let $\si=\si_o^X$.
Let $u\!:\!  X\ra \m R$ be a bounded Lipschitz superharmonic function.\\
$a)$\! For $\!\,\si$-almost all $\xi\!\in\!\bX$, the non-tangential limit $\psi(\xi)\!\!:=\!\!\underset{x\ra\xi}{\rm NTlim}\,u(x)$ exists.\\ 
$b)$ If this limit $\psi(\xi)$ is $\si$-almost surely null, then one has $u\geq 0$.
\ep

Note that the  Lipschitz continuity of $u$ 
is true for  all bounded harmonic functions, because of the Harnack inequality in Fact \ref{facHarnack1}.

\begin{proof}[Proof of Proposition \ref{pronontangentiallimit2}]~: 
This follows from the Fatou--Na\"{\i}m--Doob theorem and the Brelot-Doob trick, as they are explained  by Ancona in \cite{AnconaStFlour}.

$a)$ 
It is proven in 
\cite[Thm 1.8]{AnconaStFlour}
that for any superharmonic function $u$ on $X$ which is bounded below, 
for $\si$-almost all $\xi$ in $\bX$, the minimal fine limit $$
\psi(\xi)\!:=\!\underset{x\ra\xi}{\rm MFlim}\,u(x)
$$ 
exists.
This means that the limit of $u(x)$ when $x\ra \xi$ exists as soon as $x$ avoids a subset $E=E_\xi$ which is minimally thin at $\xi$. We recall that a subset $E\subset X$ is minimally thin if the function $P_\xi{\bf 1}_E$ is bounded by a potential on $X$. And we recall that a potential 
is a positive superharmonic function whose largest harmonic minorant is zero.
Moreover, there is a formula for this limit~:
$$
\psi(\xi)=\frac{{\rm d}\mu_h}{{\rm d}\si}(\xi)
$$
where $\mu_h$ is the trace measure on $\bX$ of the harmonic 
function $h$ in the Riesz decomposition of $u$
as a sum $u=p+h$ of a potential $p$ and a harmonic function $h$.

It is also proven in 
\cite[p.99-100]{AnconaStFlour} 
that for a Lipschitz continuous function $u$ on $X$ and a point $\xi$ in $\bX$, if the minimal fine limit $\ell\!:=\!\underset{x\ra\xi}{\rm MFlim}\,u(x)$ exists 
then the non-tangential limit 
$\!\underset{x\ra\xi}{\rm NTlim}\,u(x)$ exists 
and is equal to $\ell$.

$b)$ Since $u$ and hence its harmonic part $h$ are bounded on $X$, the measure 
$\mu_h$ is absolutely continuous to $\si$.
Hence, when the limit $\psi(\xi)$ is $\si$-almost surely zero, 
the trace measure $\mu_h$ is  zero, and the harmonic function $h$ is also $0$.
This tells us that $u$ is a potential, so that one has in particular
 $u\geq 0$.
\end{proof}

\section{The boundary transform}
\label{secboundarymap} 

\bq 
In this third chapter, we construct the boundary transform $\beta:\mc H_b( X,Y)\longrightarrow L^\infty(\bX,Y)$ and prove that it is injective.
\eq

We recall that $X$ is a \GGG Riemannian manifold, that $Y$ is a complete CAT(0) space, and that $\ol X=X\cup\partial X$.

\subsection{Harmonic maps and subharmonic functions}
\label{secharmonicsubharmonic}
\bq 
The following lemmas relate harmonic maps $h: X\ra Y$ with 
Lipschitz subharmonic functions $u$ on $ X$. 
They will allow us to apply the results on superharmonic functions
from Chapter \ref{secpreliminaryresult}.
\eq 

We begin by a useful bound for the Lipschitz constant of a harmonic map due to Cheng. 

\bl
\label{lemcheng}
Let $X$ be a Riemannian manifold with bounded geometry, and let  $Y$ be a bounded {\rm CAT(0)}-space.
There exists $L>0$ such that
for all $x_0$ in $X$, all $r\leq 1$ and any harmonic map $h:B(x_0,r)\ra Y$,
the restriction of $h$ to the ball $B(x_0,r/2)$
is $L/r$-Lipschitz.
\el

\begin{proof}
When $Y$ is a manifold,
this is a simplified version of \cite[Formula 2.9]{Cheng80}. 
See also \cite[Theorem 6]{GiaquintaHildebrandt82}.
When $Y$ is a more  general CAT(0) space, the extension of Cheng Lemma has been proven in \cite[Theorem 1.4]{ZhangZhongZhu}.
\end{proof}

\bl
\label{lemharmonicsubharmonic}
Let $ X$ be a  complete Riemannian manifold with bounded sectional curvature, and let  $Y$ be a  {\rm CAT(0)}-space.\\
$a)$ Let $h: X\ra Y$ be a bounded harmonic map and $y_0\in Y$. Then the function 
$x\mapsto d(y_0,h(x))$ is a bounded Lipschitz subharmonic function on $ X$.\\
$b)$ Let $h,h': X\ra Y$ be two bounded harmonic maps. Then the function 
$x\mapsto d(h(x),h'(x))$ is also a  bounded Lipschitz subharmonic function on $ X$.
\el

\begin{proof} When $Y$
is a manifold this is in\cite[Lemmas 3.8.1 and 3.8.2]{Jost06}.

$a)$ We can assume that the CAT$(0)$ space $Y$ is bounded.
Since $Y$ is CAT$(0)$, the function $\al$ on $Y$ defined  
by $ \al(y):=d(y_0,y)$ is convex. 
Therefore, by \cite[Lemma 1.7.1]{Jost84}
when $Y$ is a manifold
and \cite[Lemma 10.2]{EellsFuglede} in general, the function $u:=\al\circ h$ is subharmonic on $ X$.
The Lipschitz continuity of $u$ follows from the Cheng bound
in Lemma \ref{lemcheng}. 

$b)$ The proof is as in $a)$. Indeed, the map $(h,h'):  X\ra Y\times Y$ is harmonic, the product space 
$Y\times Y$ is ${\rm CAT}(0)$,
and the function $(y,y')\to d(y,y')$ is a convex function.
\end{proof}

\subsection{Construction of the boundary map}
\label{secconstructionboundary}

\bq 
In this section, we prove Proposition \ref{pronontangentiallimit}.
\eq

\begin{proof}[Proof of Proposition \ref{pronontangentiallimit}] Remember that $Y$ is here assumed to be proper. 
Fix $o\in X$ and set $\si=\si_o^X$. Let $h: X\ra Y$ be a bounded harmonic map. 
We want to prove that, for $\si$-almost all $\xi\in \bX$,
the map $h$ has a non-tangential limit $\ph(\xi)$ at the point $\xi$.

Let $y\in Y$.  By Lemma \ref{lemharmonicsubharmonic}, the function $u_y:x\ra d(y,h(x))$  is a bounded  Lipschitz subharmonic function on $ X$. 
Hence, by Proposition \ref{pronontangentiallimit2}.$a$, there exists a subset $F_y$ 
of full measure in $\bX$ such that 
the function $u_y$ admits a non-tangential limit $\psi_y(\xi)$ at each  point $\xi\in F_y$.

Let $Y_1\subset Y$ be the closure of the convex hull of $h(X)$ in $Y$. This subspace 
$Y_1$ is a bounded separable complete CAT(0) space.
Let $D\subset Y_1$ be a countable dense subset of $Y_1$. 
The intersection $F\subset\bX$ of all the sets $F_y$, for $y$ in $D$, still has full $\si$-measure. 
Note that, for all $y$, $y'$ in $Y_1$ and $x$ in $ X$, one has 
$$
d(u_y(x),u_{y'}(x))\leq d(y,y').
$$
Therefore, for all $\xi\in F$ and all $y$ in $Y_1$,
the function $u_y$ has a non-tangential limit at the point $\xi$.

We introduce the map 
\begin{eqnarray*}
\Phi:Y_1&\ra &{\rm Lip_1}(Y_1,[0,\de_{Y_1}]) \\
y'&\mapsto& (d(y,y'))_{y\in Y_{_1}}.
\end{eqnarray*}
where $\de_{Y_1}$ is the diameter of $Y_1$ and ${\rm Lip}_1$ refers to the set of $1$-Lipschitz functions endowed with the $\sup$ distance.
This map $\Phi$ is an
isometry onto its image $\Ph(Y_1)$ and, since $Y_1$ is complete, this image $\Phi(Y_1)$ is closed. 
Let $\xi\in F$. Since $Y$ is proper, the set $Y_1$ is compact, and what we have just seen tells us that the map $\Ph\circ h$
has a non-tangential limit  at the point $\xi$. Therefore,
the map $h$ also has a non-tangential limit $\ph(\xi)\in Y_1$ at the point $\xi$.
\end{proof}

\subsection{Injectivity of the boundary transform}
\label{secinjectiveboundary}
\bq 
In this section, we prove Proposition \ref{proinjectiveboundary}.
\eq 

\begin{proof} 
Let $h$ and $h'$ be two harmonic maps from $ X$ to $ Y$
whose boundary maps $\be h$ and $\be h'$ are $\si$-almost surely equal.
We want to prove that $h=h'$.

By Lemma \ref{lemharmonicsubharmonic}, the function $u:x\to d(h(x),h'(x))$
is a bounded Lipschitz subharmonic function on $ X$.
By assumption the non-tangential limit $\underset{x\ra\xi}{\rm NTlim}\,u(x)$ is zero  for $\si$-almost all $\xi$ in $\bX$.
Therefore, by Proposition \ref{pronontangentiallimit2}.$b$, the function $u$ must be non-positive. 
Since $u$ is already non-negative, we must have $u=0$, and hence $h=h'$.
\end{proof}

\section{The Poisson transform}
\label{secpoissontransform} 
\bq 
In this fourth chapter, we construct the Poisson transform.
\eq

\subsection{Density  of the Lipschitz maps}
\label{seclipschitzdense}

We first need a lemma on the density of Lipschitz maps on a compact manifold $S$ inside the set of bounded measurable maps. 
We will apply it to the boundary $S=\bOm$ 
of a bounded Lipschitz domain  $\Om$. Such a boundary is bi-Lipschitz
homeomorphic to a compact smooth manifold.

\bl
\label{lemlipschitzdense}
Let $S$ be a compact manifold and 
$Y$ be a  CAT(0) space. Then, every 
continuous map $\ph:S\ra Y$ is a uniform limit of Lipschitz maps.
\el

We begin by recalling  the classical construction of the 
weighted barycenter $\be=\be_\mu(y_0,\ldots,y_n)$ of $n\!+\!1$ points $(y_0,\ldots ,y_n)$ in a CAT(0)-space
$Y$. The weight $\mu$  belongs
to the standard $n$-simplex 
$$
\Si_n:=\{\mu=(\mu_1,\ldots,\mu_n)\mid 
\mu_i\geq 0\; \mbox{\rm for all $i$, and} \;\mu_0+\cdots+\mu_n=1\}.
$$
We endow this simplex with the $\ell^1$-distance.
This barycenter $\be$ is the unique point where the strictly convex function on $Y$ 
$$
y\mapsto \psi_\mu(y):=\sum_{0\leq i\leq n} \mu_i\,d(y_i,y)^2
$$ 
achieves its minimum. 
As a function of the weight, this barycenter map
$$
\mu\mapsto \be_\mu(y_0,\ldots,y_n)
$$ 
is $L$-Lipschitz
continuous where $L$ is the diameter 
of the finite set $\{y_0,\ldots, y_n\}$.
We refer to 
\cite[Lemma 4.2]{Kleiner99} for these properties.

\begin{proof}[Proof of Lemma \ref{lemlipschitzdense}]
Using a triangulation of $S$ we can assume that $S$
is a compact $CW$-complex. 
We endow each $n$-simplex $\Si_0$ of $S$ with the $\ell^1$-norm and we endow $S$ with the corresponding length metric. This new metric is Lipschitz equivalent to the Riemannian metric on $S$.
Each $n$-simplex $\Si_0$
of $S$ can be 
decomposed as a union of $2^n$ half-size $n$-simplices.
Iterating $k$ times this process we obtain 
a decomposition of $\Si_0$ as a union 
of $2^{kn}$ $n$-simplices of level $k$ whose size is 
$2^{-k}$ the size of $\Si_0$.

Fix $\eps>0$ and $\ph\in C(S,Y)$.  There exists an integer $k$ such that,
for each simplex $\Si$ of level $k$, one can uniformly bound the diameter
$$
{\rm diam} (\ph(\Si))\leq \eps/2.
$$
For each simplex $\Si$ of level $k$, 
we denote by $f_\Si:\Si\ra Y$ the barycenter map 
such that $f_\Si(s)=\ph(s)$ for each vertex $s$ of $\Si$.
One then has
$$
{\rm diam} (f_\Si( \Si))\leq \eps/2,
$$
and  $d(f_\Si(s),\ph(s))\leq \eps$ holds for all $s$ in $\Si$.
These maps $f_\Si$ are $2^k\eps$-Lipschitz continuous. 

These maps $f_\Si$ being compatible, each $f_\Si$ is the restriction to $\Si$ 
of a map $f:S\ra Y$. This map $f$ is also 
$2^k\eps$-Lipschitz continuous, and one has 
$d(f(s),\ph(s))\leq \eps$ for all $s$ in $S$.
\end{proof}

\bl
\label{lemcontinuousdense}
Let $S$ be a compact metric space, 
$\si$ be a Borel probability measure on $S$, and 
$Y$ be a  CAT(0) space.
Then the set $C(S,Y)$ of  continuous maps 
$f:S\ra Y$ is dense in the set $L^\infty(S,Y)$ 
of bounded measurable maps $\ph:S\ra Y$.
\el

We recall that $L^\infty(S,Y)$ is endowed with the ``topology of the convergence in probability''. The distance between two maps $\ph$, $\ph'$ in $L^\infty(S,Y)$ is given by 
\begin{equation}
	\label{eqndphphp}
	d(\ph,\ph'):=\inf\{\de\geq 0\mid\si (\{\xi\in S\mid d(\ph(\xi),\ph'(\xi))\geq \de\} )\leq \de\}\, .
\end{equation}
The space   $L^\infty(S,Y)$ and its topology
do not depend on the choice of the measure $\si$ inside its equivalence class of measures.
\vs

\begin{proof}[Proof of Lemma \ref{lemcontinuousdense}]
Let  $\ph:S\ra Y$ be a bounded measurable map. 
Let $\varepsilon >0$. By Lusin's theorem, there exists a compact subset $K\subset S$ such that the complement $K^c$ satisfies $\si(K^c)\leq \eps$ and such that the restriction $\ph|_K$ is continuous. Since a CAT(0) space $Y$ is an absolute retract metric
space, see \cite[Lemma 1.1]{Ontaneda}, and since an absolute retract metric space
is an absolute extension metric space,
there exists a continuous function $f:S\ra Y$ whose restriction to $K$
is equal to $\ph|_K$, so that $d(\varphi,f)\leq\varepsilon$.

Lusin's theorem is usually stated  under the assumption that the metric target space $Y$ is separable. Here we do not need this assumption since 
$S$ is a standard Borel space endowed with a Radon measure $\si$. Indeed, in this case, all measurable maps
$\ph: S\ra Y$ are strongly measurable. This means that $\ph$ is an almost sure limit of a sequence of measurable step functions $\ph_n$ or equivalently that there exists a conull subset $S'\subset S$ 
such that the image $\ph(S')$ is separable.	
\end{proof}

\subsection{The continuous Dirichlet problem}
\label{secdirichletproblem}
\bq
In this section we prove Proposition
\ref{prodirichlet}.
\eq

We first deal with bounded domains.

\bp
\label{prodirichleti}
Let $\Om$ be a bounded Lipschitz Riemannian domain, $Y$  a complete ${\rm CAT}(0)$ space
and $\ph: \bOm\ra Y$ a continuous map, then there exists a unique harmonic map
$h: \Om \ra Y$ which is a continuous extension of $\ph$.
\ep

\begin{proof} 
By Lemma \ref{lemlipschitzdense}, we can choose a sequence 
$\ph_n\in {\rm Lip}(\bOm,Y)$ that converges uniformly to $\ph$.
It suffices to prove that the sequence of their harmonic extensions $h_n:=P\ph_n$ given by Fact \ref{fachamilton}
converge uniformly.
We introduce the subharmonic functions
on $ \Om$ given by 
$$
u_{m,n}(x):=d(h_m(x),h_n(x)).
$$
They extend the continuous functions
on $\bOm$ given by 
$$
\psi_{m,n}(\xi)=d(\ph_m(\xi),\ph_n(\xi)).
$$
By the maximum principle, the supremum of $u_{m,n}$ on $ \Om$
is equal to the supremum of  $\psi_{m,n}$ on $\bOm$. 
Hence it goes to zero when $m,n\to\infty$.
Therefore the sequence $h_n$ converges uniformly to a map $h$
which is harmonic and which extends continuously $\ph$.

This harmonic extension is unique because if $h'$ is another harmonic extension, the positive function 
$$
x\mapsto u(x):=d(h(x),h'(x))
$$ 
is subharmonic on $\Om$ and goes to zero near the boundary. 
Hence $u=0$ and $h=h'$.
\end{proof}

We now deal with a \GGG Riemannian manifold $X$.

\begin{proof}[Proof of Proposition \ref{prodirichlet} ] 
We fix $o$ in $X$. We can assume that the diameter $\de_Y$ of $Y$ is finite. Indeed we can always replace $Y$ by a closed ball $B(o,R)$ that 
contains the image of $\ph$, since such a ball is also a 
complete ${\rm CAT}(0)$ space.

As we have seen in the proof of Lemma \ref{lemcontinuousdense},
since the compactification $\ol X$ is a metrizable compact space 
and $Y$ is a  ${\rm CAT}(0)$ space, there exists a continuous function 
$$
\psi:\ol X\ra Y
\;\;\mbox{\rm such that}\;\;
\psi|_{\bX}=\ph.
$$

We choose an increasing sequence of bounded Lipschitz domain 
$\Om_N\subset X$ such that $o\in \Om_0$ and $\Om_N$ contains the 
$1$ neighborhood of $\Om_{N-1}$.
We denote by $h_N:\Om_N\ra Y$ the harmonic extension 
of the function $\psi_N:=\psi|_{\bOm_{_N}}$
given by Proposition \ref{prodirichleti}.
We claim that 
\begin{equation}
\label{eqnclaimdirichlet} 
\forall\eps>0,\;\exists n_0>0,\;\forall N>n>n_0,\;
\forall
x\in \Om_n,\;\; d(h_N(x),h_n(x))\leq \eps \, .
\end{equation}
Since the function
$x\mapsto d(h_n(x),h_N(x))$ is subharmonic,
by the maximum principle
it is enough to check \eqref{eqnclaimdirichlet} 
for $x$ in $\bOm_n$, that is~:
\begin{equation}
\label{eqnclaimdirichlet1} 
\forall\eps>0,\exists n_0>0,\forall N>n>n_0,\;
\forall x\in \bOm_n,\;\;
d(h_N(x),\psi(x))\leq \eps \, .
\end{equation}

Let $\eps >0$. According to Lemma 
\ref{lemcontrolmeasure}, 
there  exists  
$t_0>0$ such that
for all  $x$ in $\Om_N$,
\begin{equation}
\label{eqnclaimdirichlet3} 
\si^{\Om_N}_{x}(\mc H^{t_0}_{x}(o))
\leq \eps/(2 \de_Y).
\end{equation}
By uniform continuity of $\psi$ there exists $t_1>0$ such that, 
for all $x$ in $\ol X$~: 
\begin{equation}
\label{eqnclaimdirichlet2} 
d(\psi(x),\psi(y))\leq \eps/2
\;\;\mbox{\rm for all  $y$ in $\mc H_o^{t_1}(x)$.}
\end{equation}

We choose $n_0\geq t_0\!+\!t_1$  and let $N\geq n\geq n_0$.
We fix $x$ in $\bOm_ n$
and  introduce  the subharmonic function 
$z\mapsto u(z):=d(h_N(z),\psi (x))$
on $\Om_N$.
We want to prove that $u(x)\leq \eps$.
We observe that
\begin{equation}
\label{eqnuxintuy}
u(x)\;\leq\; \int_{\bOm_N}u(y)\,
{\rm d}\si^{\Om_N}_{x}(y).
\end{equation}
Since $d(o,x)\geq  n \geq t_0+t_1$, by \eqref{eqnhsxyhtyx}, one has  
$$
X=\mc H^{t_0}_x(o)\cup \mc H_o^{t_1}(x),
$$
and we can bound this integral \eqref{eqnuxintuy}  by the sum
$I'+I''$ where~:\\ 
- $I'$ is the integral on the half-space
$\mc H^{t_0}_x(o)$,
 which by \eqref{eqnclaimdirichlet3} has small volume,\\ 
- $I''$ is the integral 
on  $\mc H^{t_1}_o(x)$ on which by \eqref{eqnclaimdirichlet2}
the function  $u$ is small.
Hence 
$$
u(x)\;\leq\;
\de_Y \;\si^{\Om_N}_{x}(\mc H^{t_0}_x(o))
\;+\;
\frac{\eps}{2}\;
\si^{\Om_N}_{x}(\mc H^{t_1}_o(x))
\;\leq\; \frac{\eps}{2}+\frac{\eps}{2}
\;=\; \eps\, ,
$$
which proves our claim \eqref{eqnclaimdirichlet}.

Now the claim \eqref{eqnclaimdirichlet}
proves that the sequence of maps $(h_N)$ converges uniformly to a harmonic map $h:X\ra Y$ that extends continuously
 $\ph$.

The proof of uniqueness is as for 
Proposition \ref{prodirichleti}.
\end{proof}

\subsection{Construction of the Poisson transform}
\label{secconstructionpoisson}

The construction uses the following classical ``continuous extension theorem''.

\bl 
\label{lemuniconext}
Let $E$ be a metric space, $D\subset E$ a dense subset and $F$ a complete metric space.
Then every uniformly continuous map $P:D\ra F$ admits a unique continuous extension $P:E\ra F$. 
\el
\begin{proof}[Proof of Lemma \ref{lemuniconext}] 
	This is classical.
\end{proof}

\begin{proof}[Proof of Proposition \ref{propoissontransform}]
Remember that $Y$ is assumed to be bounded.
We use Lemma \ref{lemuniconext} 
with $E=L^\infty(\bX,Y)$, $D=C(\bX,Y)$
and $F=\mc H_b( X,Y)$. Note that $F$ is a complete 
metric space since a uniform limit of harmonic maps is harmonic.

We want to prove that the map $P:C(\bX,Y)\ra \mc H_b( X,Y)$ given by Proposition \ref{prodirichlet}
has a unique continuous extension to $L^\infty(\bX,Y)$.
By Lemmas \ref{lemuniconext} and \ref{lemcontinuousdense}, it suffices to prove that this map 
$P$ is uniformly continuous.
We fix a compact $K\subset X$ and a point $o\in K$, and we set 
$$
C_K=\sup\limits_{\xi\in \bX,\, x\in K}\! P_\xi(x)
\;<\; \infty\, ,
$$ 
where 
$P_\xi(x)$ is the Poisson kernel.

Let $0<\eps\leq 1$ and $\ph$, $\ph'$ be two  continuous maps from $\bX$ to $Y$ 
such that $d(\ph,\ph')\leq \eps$. 
This means that the function $\psi$ on $\bX$
defined, for $\xi$ in $\bX$, by $\psi(\xi):=d(\ph(\xi),\ph'(\xi))$ satisfies
\begin{equation} 
\label{eqndistanceproba}
\si_{o}(\{\xi\in \bX\mid \psi(\xi)\geq \eps\})\leq \eps,
\;\;\mbox{\rm where $\si_o=\si_o^X$}.
\end{equation}
Note  that this function $\psi$ is bounded by the diameter $\de_Y$ of $Y$.

Let $h=P\ph$ and $h'=P\ph'$ 
be their harmonic extensions to $X$. 

By Lemma \ref{lemharmonicsubharmonic}, the continuous function $u$ on $X$ 
given, for $x$ in $X$,  by 
$u(x)\!:=\! d(h(x),h'(x))$ is  subharmonic on $ X$. This function is a continuous extension of $\psi$.
Therefore it satisfies, for all $x$ in $X$,
\begin{eqnarray*}
u(x)&\leq &\int_{\bX}  \psi(\xi)\,P_\xi(x)\,{\rm d}\si_o(\xi)\, .
\end{eqnarray*}
Plugging \eqref{eqndistanceproba} in this inequality, one gets for every $x$ in $K$~:
\begin{eqnarray*}
u(x)&\leq &\eps\int_{\bX} P_\xi(x)\,{\rm d}\si_o(\xi) +\de_Y\int_{\{\psi(\xi)\geq \eps\}} P_\xi(x)\,{\rm d}\si_o(\xi)\\
&\leq& \eps + C_K\de_Y\eps.
\end{eqnarray*}
This proves, for any $\varphi,\varphi'$ in $C(\partial X,Y)$ and any compact subset $K\subset X$, the inequality~:
$$
\sup_{x\in K}d(h(x),h'(x))\leq (1+C_K\de_Y)d(\ph,\ph')
\, .
$$
This is the uniform continuity of the map $P$.
\end{proof}

\section{The boundary and the Poisson transform}
\label{secboundarypoisson}

\bq
In this chapter, we prove that the Poisson transform $P$ is the inverse of the boundary transform $\be$.
\eq 

\subsection{The Lebesgue density theorem}
\label{seclebesgue}

We recall here the generalized Lebesgue density theorem. 
See \cite[Section 4.6]{BenoistFive} for a short and complete proof. 
\bd
\label{defbdistance}
A quasi-distance on a space $S$ is a map 
$d_0 : S\times S \ra [0, \infty[$ for which there exists 
$b>0$ \\ 
such that
$
d_0(\xi_1, \xi_3) \leq 
b (d_0(\xi_1, \xi_2)+d_0(\xi_2, \xi_3))\; ,\; \;\forall  \xi_1,\xi_2,\xi_3 \in S$,\\
such that $d_0(\xi_1, \xi_2) = d_0(\xi_2, \xi_1)$ and\\
such that
$d_0(\xi_1, \xi_2) = 0 \Leftrightarrow \xi_1 = \xi_2$.
\ed

- In this case, one says that $S$ is a quasi-metric space. Then, there exists  a topology on
$S$ for which the balls 
$B(\xi, \eps) := \{ \eta \in S \mid d_0(\xi, \eta) \leq \eps \}$,
with $\xi\in S$ and $\eps > 0$, form a basis of
neighborhood of the points $\xi$.\\
- One then has the inclusion $\ol{B(\xi, \eps)} \subset B(\xi, b \eps)$.
\vs 

Let $(S, d_0)$ be a  compact quasi-metric space and 
$\si$ be a finite Borel measure on $S$. One says that 
$\si$ is doubling 
if there exists  $C>0$ such that, for all $\xi\in S$ and 
$r > 0$, one has $\si(B(\xi, 2r)) \leq C\,\si(B(\xi,r))$.

Let $F\subset S$ be a measurable subset. A point $\xi \in S$ is called a density point if
$$\lim\limits_{\eps\ra 0}\frac{\si (B(\xi, \eps) \cap F)}{\si(
B(\xi, \eps))}
= 1.
$$

\bfa
{\bf (Lebesgue)} 
\label{faclebesgue}
Let $(S, d_0)$ be a compact quasi-metric space,  $\si$ be a doubling finite Borel measure on $S$, and let $F$ be a measurable subset of $S$.
Then $\si$-almost every point of $F$ is a density point.
\efa

In the next section, we will apply Fact \ref{faclebesgue} with 
$S=\bX$ and with $\si=\si_{o}$.

We will use the  quasi-distance on $\bX$  defined, for two points 
$\eta_1$ and $\eta_2$, by 
the exponential inverse of the Gromov product~:
\begin{equation}
\label{eqnbdistance}
d_0(\eta_1,\eta_2)= e^{-(\eta_1|\eta_2)_{o}}.
\end{equation}
This formula defines indeed a quasi-distance on $\bX$, 
because of \eqref{eqngromovproduct}.
Note that one can modify this formula so that $d_0$ is actually a distance, 
see \cite{GhysHarp90}.

The balls for this quasi-distance in $\bX$ are 
the 
trace at infinity of the 
half-spaces of $X$. 
The doubling property for 
$\si_{o}$ is proven in Lemma \ref{lemsidoubling}.

\subsection{Limit of subharmonic functions}
\label{seclimitsubharmonic}

\bq 
In this section we prove the technical lemma \ref{lemradialconvergence}
that plays a crucial role in the proof of 
the surjectivity of the boundary transform.
\eq 

We fix a point $o$ in $X$. For all $\xi\in \bX$ we define $N\xi$ as the union 
of the   geodesic rays $o\xi$
from $o$ to $\xi$.
We define then  {\it the  tube $NF$ 
over a compact set} $F\subset \bX$ as 
the union 
\begin{eqnarray} 
\label{eqntube} 
NF:=\bigcup_{ \xi\in F}N\xi \subset  X.
\end{eqnarray}

\bl 
\label{lemradialconvergence} 
Let $ X$ be a \GGG Riemannian manifold, and fix a point $o\in X$.
Let  $\psi_n\!:\!\bX\!\ra\! [0,1]$ be a sequence of con\-ti\-nuous functions that converges $\si_{o}$-almost surely to $0$.\!\!
Let $u_n\!:\!X\!\ra\! [0,1]$ be non-negative subharmonic functions on $X$ 
that extend continuously $\psi_n$. 
Then, for all $\eps>0$, there exists a compact subset 
$F_{\eps}\subset\bX$ with
$\si_{o}( F_{\eps}^c)\leq\eps$,
such that the sequence $(u_n)$ converges uniformly to $0$ on 
the  tube $NF_\eps$.
\el

The arguments of Section \ref{secconstructionpoisson} tell us that 
the sequence $(u_n)$ converges to $0$ uniformly on the compact subsets $K$
of $ X$. Lemma \ref{lemradialconvergence} tells us that this convergence 
is still uniform on ``large radial subsets of $ X$''.

\begin{proof} 
{\bf First step } {\it We control the  Poisson kernel on tubes.}
 
For $\xi$ in $\bX$ and  $m\geq 0$ we 
denote by 
$$
B_m(\xi):=B(\xi,e^{-m})=
\{\eta\in \bX\mid (\eta|\xi)_o\geq m\}
$$ 
the balls for the 
quasidistance $d_0$  in \eqref{eqnbdistance},
and we introduce the annuli
$$
A_m(\xi):=B_m(\xi)\smallsetminus B_{m+1}(\xi),
$$
so that one has
$\cup_{m\geq 0}A_m(\xi)=\bX\!\smallsetminus\!  \{\xi\}$.

By Proposition \ref{propoisson2}, there exists  $C_3>0$ such that,
for  $\xi\in \bX$ and $m\geq 0$, 
\begin{equation} 
\label{eqnratiopoisson2}
\frac{P_{\eta_{_1}}(x)}{P_{\eta_{_2}}(x)}
\;\leq\; C_3
\;\;\mbox{\rm holds for all $x\in N\xi$ and $\eta_1$, $\eta_2$ in $A_m(\xi)$.} 
\end{equation}

{\bf Second step} {\it We apply the Lebesgue density theorem.}

Let $\eps>0$. Since the sequence $(\psi_n)$ converges $\si_{o}$-almost surely to $0$, 
there exist an integer $n_\eps\geq 1$ and a compact subset $K_\eps\subset \bX$ with
$\si_{o}(K_\eps^c)\leq \eps/2$, and such that 
 $\psi_n(\xi)\leq \eps$ for all $n\geq n_\eps$ and $\xi\in K_\eps$.

By the Lebesgue density theorem (Fact \ref{faclebesgue}), applied to the doubling measures $\si_{o}=\si_o^X$
and the family of balls $B_m(\xi)$,
the sequence of functions 
$$
f_m(\xi):=\frac{\si_{o}(B_m(\xi)\cap K_\eps^c )}{ \si_{o}(B_m(\xi))}
$$
converges to zero for $\si_{o}$-almost all $\xi\in K_\eps$. 
Since the measure $\si_{o}$ is doubling, the ratios 
$\displaystyle
\frac{\si_{o}( B_m(\xi))}{\si_{o}( B_{m+1}(\xi))}
$ 
are uniformly bounded,
and this implies that the ratios 
$\displaystyle 
\frac{\si_{o}( B_m(\xi))}{\si_{o}(A_m(\xi) )}$  
are also uniformly bounded.
Hence the sequence of functions
$\displaystyle
g_m(\xi):=\frac{\si_{o}(A_m(\xi)\cap K_\eps^c)}{ \si_{o}(A_m(\xi))}
$
also converges to zero for $\si_{o}$-almost all $\xi\in K_\eps$.

Therefore, by Egorov theorem, there exist a compact subset $L_\eps\subset K_\eps$ 
and an integer $m_\eps\geq 1$ such that
$\si_{o}(L_\eps^c)\leq \eps$  and with
\begin{equation} 
\label{eqnratiomassanuli}
\si_{o}(A_m(\xi)\cap K_\eps^c)\leq \eps\, \si_{o}(A_m(\xi))\, 
\;\;\mbox{\rm for all $m\geq m_\eps$ and $\xi\in L_\eps$.} 
\end{equation}

{\bf Third step} {\it We bound the functions $u_n$ by using  the Poisson kernel.}

Since each function $u_n$ is subharmonic with boundary value $\psi_n$, one has 
\begin{eqnarray*}
u_n(x)&\leq &\int_{\bX} \psi_n(\eta)\,P_\eta(x) \,{\rm d}\si_{o}(\eta)
\;\;\mbox{\rm for all $x\in  X$.} 
\end{eqnarray*}
We now assume  that $x$ belongs to a tube $N\xi$ 
with $\xi\in L_\eps$. We write
\begin{eqnarray*}
u_n(x)&\leq &\sum_{m=0}^{\infty}I_{m,n}(x,\xi)
\;\; \mbox{\rm where 
$I_{m,n}(x,\xi):=\displaystyle 
\int_{A_m(\xi)} \psi_n(\eta)\,P_\eta(x) \,{\rm d}\si_{o}(\eta)$.} 
\end{eqnarray*}
We split this sum into two parts, according to whether $m< m_\eps$ or $m\geq m_\eps$.

First assume that $m<m_\eps$.
The function $(x,\eta)\mapsto P_\eta(\xi)$ 
being continuous on $(\ol{X}\times\bX)\smallsetminus \De_{\bX}$, there exists a constant $C_8=C_8(m_\eps)>0$ 
such that one has, for all $\xi\in\bX$~:
\begin{eqnarray*}
\label{eqnpxeta}
P_{\eta}(x)\leq C_8
\;\;\;\;\mbox{\rm  for all $x\in N\xi$
and $\eta\in \bX\smallsetminus B_{m_{_\eps}}(\xi)$.}
\end{eqnarray*}
This gives the bound
\begin{eqnarray}
\label{eqnksmall}
\sum_{m<m_\eps}I_{m,n}(x,\xi)&\leq& C_8
\int_{\bX} \psi_n(\eta)\,{\rm d}\si_{o}(\eta)\, ,
\end{eqnarray}
where the integral converges to $0$ when $n\to\infty$ by the Lebesgue dominated convergence theorem.

Now assume that $m\geq m_\eps$. One splits the integral $I_{m,n}(x,\xi)$ as a sum
\begin{eqnarray*}
\label{eqnipkn0}
I_{m,n}(x,\xi)&=&I'_{m,n}(x,\xi)+I''_{m,n}(x,\xi)
\;\; \mbox{\rm where}\\
I'_{m,n}(x,\xi)&:=& 
\int_{A_m(\xi)\cap K_\eps^c} \psi_n(\eta)\,P_\eta(x) \,
{\rm d}\si_{o}(\eta),\\
I''_{m,n}(x,\xi)&:=& 
\int_{A_m(\xi)\cap K_\eps} \psi_n(\eta)\,P_\eta(x) \,
{\rm d}\si_{o}(\eta).
\end{eqnarray*}
Since $m\geq m_\eps$ and $\xi\in L_\eps$, we obtain by using   \eqref{eqnratiopoisson2} and \eqref{eqnratiomassanuli} and the bound $\|\psi_n\|_\infty \leq 1$~:
\begin{eqnarray}
\label{eqnipkn}
\nonumber 
I'_{m,n}(x,\xi)&\leq& \max_{\eta\in A_m(\xi)}P_\eta(x)\; \; 
\si_{o}(A_m(\xi)\cap K_\eps^c)\\
\nonumber
&\leq& C_3\;\min_{\eta\in A_m(\xi)}P_\eta(x)\;\; \eps\,\si_{o}(A_m(\xi))\\
&\leq& 
\eps \, C_3\int_{A_m(\xi)} P_\eta(x) \,{\rm d}\si_{o}(\eta).
\end{eqnarray}
Assume moreover that $n\geq n_\eps$. Using the definition of $K_\eps$, we obtain~: 
\begin{eqnarray}
\label{eqnippkn}
I''_{m,n}(x,\xi)
&\leq& 
\eps \int_{A_m(\xi)} P_\eta(x) \,{\rm d}\si_{o}(\eta).
\end{eqnarray}
Combining \eqref{eqnipkn} and \eqref{eqnippkn} and summing over $m\geq m_\eps$, one gets 
for $n\geq n_\eps$~:
\begin{equation}
\label{eqnklarge}
\sum_{m\geq m_\eps}I_{m,n}(x,\xi)
\;\;\leq\;\; 
(1+C_3)\,\eps
\int_{\bX} P_\eta(x)\,{\rm d}\si_{o}(\eta)\\
\;\; =\;\;
(1+C_3)\,\eps \, .
\end{equation}
We now define the compact $F_\eps$ as the intersection 
$F_\eps:=\cap_{\ell\geq 1}L_{\eps_\ell}$ with $\eps_\ell:=2^{-\ell}\eps$, so that $\si_{o}(F_\eps^c)\leq \eps$.
Combining \eqref{eqnksmall} and \eqref{eqnklarge} we observe that one has, for all $x$ in $N F_\eps$ and every integers $\ell\geq 1$ and $n\geq n_{\eps_{_\ell}}$~:
\begin{eqnarray*} 
u_n(x)&\leq& C_8
\int_{\bX} \psi_n(\eta)\,{\rm d}\si_{o}(\eta)
\; + (1+C_3)\, \eps_\ell\, .
\end{eqnarray*}
If $\ell$ is large enough the second term is small. 
And, as we have already seen by using the Lebesgue dominated convergence theorem, the
first term is small if $n$ is large enough.

This proves that the sequence $(u_n)$ converges uniformly to $0$
on $NF_\eps$.
\end{proof}

\subsection{Surjectivity of the boundary transform}
\label{secsurjectiveboundary}

\begin{proof}[Proof of Proposition \ref{prosurjectiveboundary}]
Let $\ph\in L^\infty(\bX,Y)$. We want to prove the equality $\ph=\be P\ph$. The metric space $Y$ has been assumed to be 
proper so that we can use the existence of the boundary map $\be$ from Proposition \ref{pronontangentiallimit}. 
Let $\de_Y$ be the diameter of $Y$.

Let $(\ph_n)$ be a sequence in $ C^0(\bX,Y)$ that converges almost surely to $\ph$.
Such a sequence also converges to $\ph$ in probability i.e. for the distance \eqref{eqndphphp}.
Let $h_n=P\ph_n$ and $h=P\ph$. By construction, one has $\ph_n=\be h_n$
and the harmonic map $h$ is the limit of the harmonic maps $h_n$, where
 the convergence is uniform on compact sets of $ X$. 
This is not enough to conclude. But we will prove below that,
for all $\eps>0$,
this convergence is also uniform on the  tube $NF_\eps$ 
over a compact set 
$F_\eps\subset \bX$ such that $\si_{o}(  F_\eps^c)\leq \eps$.

Since the sequence $(\ph_n)$ converges almost surely to $\ph$, the continuous functions 
$
\psi_{m,n}:\bX\ra [0,\de_Y]$ defined, for $\xi$ in $\bX$, by
$$
\psi_{m,n}(\xi)=d(\ph_m(\xi),\ph_n(\xi))
$$ 
converge almost surely to $0$ 
when $m,n$ go to $\infty$.

The functions
$u_{m,n}:X\ra [0,\de_Y]$ defined by
$$
u_{m,n}(x)=d(h_m(x),h_n(x))
$$ 
for $x$ in $X$
extend continuously the functions $\psi_{m,n}$.
By Lemma \ref{lemharmonicsubharmonic}, these functions  $u_{m,n}$
are subharmonic on $ X$.

Let $\eps\!>\!0$. Lemma \ref{lemradialconvergence} ensures that
there exists a compact subset 
$F_{\eps}\!\subset\!\bX$ such that
$\si_{o}( F_{\eps}^c)\leq\eps$, and such that
the sequence $(u_{m,n})$ converges uniformly to $0$ on 
the  tube $NF_\eps$.
This tells us that the convergence of the sequence $(h_n)$ 
to $h$ is uniform on the tube $NF_\eps$.
By Egorov theorem, we may also assume that 
the sequence of continuous functions  $(\ph_n)$ converges uniformly to $\ph$ on $F_\eps$.
Therefore the function $\ol h:\ol X\ra Y$ equal to $h$ on $X$
and equal to $\ph$ on $\bX$ is continuous on 
$NF_\eps\cup F_\eps$.

This proves that, when $\xi$ is in $F_\eps$, 
the non-tangential limit 
$
\underset{x\ra\xi}{\rm NTlim}\,h(x)$
given by Proposition \ref{pronontangentiallimit}
is equal to $\ph(\xi)$.
Since the measure of $F_\eps^c$ is arbitrarily small,
the map $\ph$ is the boundary map of $h$.
\end{proof}

\subsection{A concrete example}

We give an example of 
Theorem \ref{thmbijectiveboundary}  in a situation 
where $X\!=\!Y$ is the hyperbolic plane $\m H^2$ and the boundary map $\ph:\partial\m H^2\ra\m H^2$ has finite image.

\bc
\label{corexample}
Let $X:=\m H^2$ and $x_1,\ldots,x_n$  be $n$  points on $\partial \m H^2$ 
cutting $\partial \m H^2$ into $n$ open arcs $I_1,\ldots, I_n$.
Let $y_1,\ldots,y_n$ be $n$ points on $Y:=\m H^2$. 
Then there exists a unique bounded harmonic map $h:\m H^2\ra \m H^2$ that extends continuously to the arcs $I_j$
with
$h(I_j)=\{y_j\}$, for all $j$.
\ec

\begin{proof}
By Theorem \ref{thmbijectiveboundary} the  map $h$ has to be the Poisson transform of the map $\ph:\partial\m H^2 \ra \m H^2$
that sends the sides $I_j$ to the point $y_j$, for all $j$.
We only have to check that this map $h$ extends continuously $\ph$ 
outside the points $x_j$. This is the content of Proposition \ref{procontinuityboundary}.
\end{proof}

\bp 
\label{procontinuityboundary}
Let $X$ be a  \GGG Riemannian manifold,\! and $Y$ be a proper {\rm CAT(0)}-space.
Let  $h\!:\! X\!\ra\! Y$ be a bounded harmonic map
and $\ph:\bX\ra Y$ be its boundary map. 
Let $I\subset \bX$ be an open set on which $\ph$ is continuous.
Then, for all $\xi\!\in\! I$, one has $\ph(\xi):=\underset{x\ra\xi}{\rm lim}\,h(x)$.
\ep

\begin{proof}
The argument is as in the proof 
of Proposition \ref{prodirichlet} in Section \ref{secdirichletproblem}.
\end{proof}

\subsection{Bounded Lipschitz domain}
\label{seclipschitzdomain}

A ``bounded Lipschitz Riemannian domain'' $ \Om$ means 
a connected bounded open subset of a smooth Riemannian manifold $(M,g_0)$
such that $ \Om$ is the interior of a connected 
compact submanifold $\ol \Om$ of $M$ whose boundary $\bOm$ is a non-empty Lipschitz continuous codimension one submanifold.
\vs 

The same argument as for Theorem \ref{thmbijectiveboundary} 
will give the following corollary

\bc
\label{corbijectiveboundary}
Let $\Om$ be  a bounded Lipschitz Riemannian domain, and $Y$ be a proper {\rm CAT(0)}-space.\\
$a)$ Let  $h\!:\! \Om\!\ra\! Y$ be a bounded harmonic map. 
Then, for $\si$-almost all $\xi\!\in\! \bOm$,
the map $h$ admits a non-tangential limit $\ph(\xi):=\underset{x\ra\xi}{\rm NTlim}\,h(x)$ at the point $\xi$.\\
$b)$ The map $h\mapsto \be(h):=\ph$
gives a bijection 
$
\beta:\mc H_b( \Om,Y)\stackrel{\sim}{\longrightarrow} L^\infty(\bOm,Y).
$
\ec
In this case, a non-tangential limit means a limit along all sequences $x_n$ such that
$\displaystyle\sup\limits_{n\geq 1}\frac{d(x_n, \xi)}{d(x_n,\bOm)}<\infty$.

The measure $\si$ is any finite Borel measure on $\bOm$ which is equivalent to any of
the harmonic measures of $\Om$.
Since $\Om$ is a bounded Lipschitz Riemannian domain, by Dahlberg's theorem in \cite{Dahlberg79}, one can 
choose $\si$ to be the Riemannian measure on $\bOm$.
\vs 

As in Section \ref{secprevious},
when $Y$ is a Riemannian manifold, Corollary \ref{corbijectiveboundary}.$a$ is due to
Aviles, Choi, Micallef in \cite[Thm 5.1]{AvilesChoiMicallef91},
and Theorem \ref{corbijectiveboundary}.$b$ is expected to be true 
as a final observation 
in \cite[Section 1]{AvilesChoiMicallef91}. 
The first cases of Corollary \ref{corbijectiveboundary}.$b$ 
that seem to be new is when 
$X$ is the Euclidean unit ball in $\m R^k$ and 
$Y$ is the hyperbolic space $\m H^\ell$.

\begin{proof}
Corollary \ref{corbijectiveboundary} 
is a corollary of the proof of 
Theorem \ref{thmbijectiveboundary}. 
The strategy is the same, relying on variations of 
Propositions \ref{propoisson2} and \ref{pronontangentiallimit2}
for bounded Lipschitz Riemannian domains. 
The proofs of these variations are very similar. The only difference is that they rely on \cite{Ancona78} instead of \cite{AnconaStFlour}. 
\end{proof}

\br
\label{remtrickBHK}
The fact that Corollary \ref{corbijectiveboundary} is a special case of Theorem \ref{thmbijectiveboundary} can also be explained 
thanks to a trick due to Bonk, Heinonen and Koskela 
in \cite[Chapter 8]{BonkHeinonenKoskela}. This trick 
consists in replacing the Riemannian metric $g_0$ on $\Om$ by $g=\ol{d}(x)^{-2}g_0$ where $\ol d$ is a suitable $C^\infty$ function roughly equal to the  distance to the boundary,
obtaining this way a \GGG Riemannian manifold $(\Om,g)$.
One then sees the harmonic and subharmonic functions on 
$(\Om,g_0)$ 
as $\mc L$-harmonic and $\mc L$-subharmonic functions 
on  $X$ where $\mc L:= \ol{d}(x)^2\De_{g_0}$ is an elliptic
differential  operator  of order $2$ which is equal to the Laplacian
$\De_g$ up to terms of order $1$ and which also has spectral gap. 
All the arguments we developed in this paper for the Laplace operator 
$\De_g$ also apply to the operator $\mc L$.
This trick could be applied to a much wider class of bounded Riemannian domains called ``inner uniform domains''. See also  Aikawa in \cite[Theorem 1.2]{Aikawa04}
for more on these domains. 
\er

{\small
	\bibliography{harmonicbounded}
}
\vs 

{\small\noindent 
	Y. Benoist:   CNRS \& Université Paris-Saclay 
	\hfill{\tt yves.benoist@u-psud.fr}\\
	D. Hulin:\;\;\;	Université Paris-Saclay 
	\hfill{\tt dominique.hulin@u-psud.fr}}

\end{document}